\documentclass[a4paper,11pt,twoside,reqno]{amsart}

\usepackage[utf8]{inputenc}
\usepackage[plainpages=false,pdfpagelabels=true]{hyperref}
\usepackage{xcolor,amssymb,amsthm,wasysym}
\usepackage{slashed}

\numberwithin{equation}{section}

\newtheorem{Satz}{Theorem}[section]
\newtheorem{Prop}[Satz]{Proposition}
\newtheorem{Lem}[Satz]{Lemma}

\newtheorem{Cor}[Satz]{Corollary}
\newcommand{\cR}{{\mathcal R}}

\theoremstyle{remark} 
\newtheorem{Bem}[Satz]{Remark}
\newtheorem{Bsp}[Satz]{Example}

\allowdisplaybreaks[1]

\renewcommand{\epsilon}{\varepsilon}

\newcommand{\R}{\ensuremath{\mathbb{R}}}

\newcommand{\Z}{\ensuremath{\mathbb{Z}}}
\newcommand{\D}{\slashed{D}}
\newcommand{\p}{\slashed{\partial}}
\newcommand{\sff}{\mathrm{I\!I}}

\title{The evolution equations for regularized Dirac-geodesics}
\author{Volker Branding}
\date{\today}
\address{TU Wien\\
Institut für diskrete Mathematik und Geometrie\\
Karlsplatz 13, 1040 Wien}
\email[]{volker@geometrie.tuwien.ac.at}
\urladdr{http://www.geometrie.tuwien.ac.at/branding/}
\begin{document}

\begin{abstract}
We study the evolution equations for a regularized version of Dirac-geodesics,
which are the one-dimensional version of Dirac-harmonic maps. We show that for the regularization
being sufficiently large, the evolution equations subconverge to a regularized Dirac-geodesic.
In the end, we discuss the limiting process of removing the regularization.
\end{abstract} 

\maketitle

\section{Introduction and Results}
Harmonic maps between Riemannian manifolds \cite{MR0164306} are well studied objects in differential geometry.
They form a nice variational problem with a rich structure. In its simplest version, harmonic maps manifest as geodesics.
The existence of closed geodesics in a closed Riemannian manifold can be studied by various methods and is well understood at present. 
A powerful method in this case is the so called heat flow method.
Here, one deforms a given curve by a heat-type equation. For geodesics, this method was successfully applied by Ottarsson \cite{MR834094}. 
The same method can also be applied when the domain manifold is not one-dimensional. In this case Eells and Sampson established their famous existence result
for harmonic maps under the assumption that the target manifold has non-positive curvature.
\par\medskip
An extension of harmonic maps introduced recently are Dirac-harmonic maps \cite{MR2262709}. 
A Dirac-harmonic map is a pair \((\phi,\psi)\), (with \(\phi\colon M\to N\) being a map and \(\psi\) a vector spinor) which
is a critical point of an energy functional that is motivated from supersymmetric field theories in physics.
In physics, this functional is known as supersymmetric sigma model, for the physical background, see \cite{MR1701598}.
In contrast to the physical literature the spinor fields we are considering here are real-valued and commuting.
The Euler-Lagrange equations for Dirac-harmonic maps couple the equation for harmonic maps with spinor fields.
They form a weak elliptic system that couples a first and a second order equation.
As limiting cases Dirac-harmonic maps contain both harmonic maps and harmonic spinors.

At present many essential results for Dirac-harmonic maps have been obtained.
The regularity theory for Dirac-harmonic maps is fully developed, see \cite{MR2544729}.
In addition, important concepts like an energy identity have also been established \cite{MR2267756}.
A classification result for Dirac-harmonic maps between surfaces was obtained in \cite{MR2496649}.

However, the existence question for Dirac-harmonic maps is still not fully understood.
Some explicit solutions of the Euler-Lagrange equations for Dirac-harmonic maps are constructed in \cite{MR2569270}, see also \cite{NicolasGinoux}.
A general existence result could be obtained in the situation that the Euler-Lagrange equations for Dirac-harmonic maps decouple \cite{ammannginoux}.
Here, the authors assume the existence of a harmonic map \(\phi_0\) and use methods from index theory to construct a spinor \(\psi\)
such that one gets a Dirac-harmonic map \((\phi_0,\psi)\). An existence result for the boundary value problem for Dirac-harmonic maps was recently obtained in \cite{springerlink:10.1007/s00526-012-0512-5}.
\par\medskip
In this article we want to treat the existence question in the most simple situation.
Namely, we assume that \(M=S^1\) and \(N\) being a closed Riemannian manifold.
Our approach uses the heat-flow method, but we cannot apply it directly since the Dirac operator is
of first order and thus we would not get a reasonable evolution equation. Hence, we first of all consider
a regularization of the energy functional for Dirac-geodesics and consider the heat-flow of the regularized functional.
Finally, we study the removal of the regularization. The basic idea of our approach is to
solve a problem that is better to deal with than Dirac-geodesics and still get an existence result for Dirac-geodesics.

A similar approach was performed in \cite{springerlink:10.1007/s00526-011-0404-0}.
There, an extra term \(F(\gamma,\psi)\) in the energy functional is considered and depending on the properties of this term
existence results are obtained.

Here, we study the energy functional
\begin{equation}
E_\epsilon(\gamma,\psi)=\frac{1}{2}\int_{S^1}(|\gamma'|^2+\langle\psi,\D\psi\rangle+\epsilon|\D\psi|^2)ds,
\end{equation}
where \(\gamma\colon S^1\to N\) and differentiation with respect to the curve parameter \(s\)
is abbreviated by \('\). Moreover, \(\psi\) is a vector spinor and \(\D\) the Dirac operator acting on \(\psi\).
This functional can be interpreted as a regularized version of the energy for Dirac-geodesics, which can be 
obtained in the limit \(\epsilon\to 0\).  The idea is to obtain an existence result for critical points of \(E_\epsilon(\gamma,\psi)\) 
and then study the limit \(\epsilon\to 0\). Of course, we want to keep the regularizing parameter \(\epsilon\) small
such that we can think of \(E_\epsilon(\gamma,\psi)\) as a small perturbation of the energy functional for Dirac-geodesics.

The Euler-Lagrange equations of \(E_\epsilon(\gamma,\psi)\) are given by:
\begin{align}
\tau(\gamma)=&\cR(\gamma,\psi)+\epsilon\cR_c(\gamma,\psi), \\
\nonumber\epsilon\tilde{\Delta}\psi=&\D\psi,
\end{align}
where \(\tau\) denotes the tension field of the curve \(\gamma\) and \(\tilde{\Delta}\) the Laplacian acting on vector spinors.
Note that \(\D^2=-\tilde{\Delta}\) since our domain is one-dimensional.
The curvature terms \(\cR(\gamma,\psi)\) and \(\cR_c(\gamma,\psi)\) are given by
\begin{align}
\label{curvature-term1}\cR(\gamma,\psi)=&\frac{1}{2}R^N(\partial_s\cdot\psi,\psi)\gamma',\\
\label{curvature-term2}\cR_c(\gamma,\psi)=& R^N(\tilde{\nabla}\psi,\psi)\gamma'.
\end{align}
Here, \(R^N\) denotes the Riemann curvature tensor on \(N\), the \(\cdot\) refers to Clifford multiplication
and \(\tilde{\nabla}\) represents the covariant derivative acting on vector spinors.

In order to obtain an existence result for critical points of \(E_\epsilon(\gamma,\psi)\) we use the \(L^2\)-gradient flow
of this functional, which is given by the following set of coupled evolution equations:
\begin{align}
\label{evolution-gamma-intro}
\frac{\partial\gamma_t}{\partial t}=&\tau(\gamma_t)-\cR(\gamma_t,\psi_t)-\epsilon\cR_c(\gamma_t,\psi_t),\\
\label{evolution-psi-intro}
\frac{\tilde{\nabla}\psi_t}{\partial t}=&\epsilon\tilde{\Delta}\psi_t-\D\psi_t
\end{align}
together with the initial data \((\gamma(s,0),\psi(s,0))=(\gamma_0(s),\psi_0(s))\).
\par
Equations (\ref{evolution-gamma-intro}) and (\ref{evolution-psi-intro}) form a 
parabolic system that behaves nicely from an analytical point of view. 
Since we take the domain manifold to be \(S^1\) we can use the Bochner technique
and the maximum principle to derive energy estimates.
\par
Finally, we will prove the following
\begin{Satz}
\label{theorem-curve}
Assume that \(M=S^1\) with fixed spin structure and \(N\)
is a compact Riemannian manifold without boundary. Then for any sufficiently regular initial data \((\gamma_0,\psi_0)\)
and \(\epsilon\) small, there exists a unique smooth solution of (\ref{evolution-gamma-intro}) and (\ref{evolution-psi-intro}) 
for all \(t\in [0,\infty)\). \\
If \(\epsilon\geq 1\), the evolution equations subconverge at infinity, that is, there exists a sequence \(t_k\to\infty\)
such that \((\gamma_{t_k},\psi_{t_k})\to(\gamma_{\infty},\psi_{\infty})\) in \(C^2\), where \((\gamma_{\infty},\psi_{\infty})\)
is a regularized Dirac-geodesic homotopic to \((\gamma_0,\psi_0)\).
\end{Satz}
In principle, we can then take the limit \(\epsilon\to 0\) and at least we do not lose any regularity.
However, by our approach it does not seem possible to further characterize the solution \((\gamma_\infty,\psi_\infty)\) after \(\epsilon\to 0\).
\par\medskip
The results presented in this article are part of the author's PhD thesis \cite{phd}.
\par\medskip
Let us now describe the framework for Dirac-geodesics in more detail.
Assume that \(M\) is a closed curve, for simplicity take \(M=S^1\). In addition, let $(N,g_{ij})$ be a compact, smooth Riemannian manifold without boundary.
We will denote the curve parameter with \(s\) and in order to shorten the notation we write \('=\frac{\partial}{\partial s}\).
Coordinates on the target manifold $N$ are denoted by $y$. We use the Einstein summation convention, which means that we will sum over repeated indices.\\

On \(S^1\) there are two spin structures, which will be abbreviated by $\sigma_1$ and $\sigma_2$.
In the case of $\sigma_1$, spinors can be identified as periodic complex-valued functions on \(S^1\)
satisfying $\psi(s+2\pi)=\psi(s)\). Regarding the other spin structure $\sigma_2$,
spinors can be identified as antiperiodic complex-valued functions on \(S^1\) satisfying $\psi(s+2\pi)=-\psi(s)$.
It is well known that Clifford multiplication on $S^1$ is given by multiplication with the imaginary unit $i$,
but in order to keep track of its skew-symmetry we write \(\partial_s\cdot\). Clifford multiplication is skew-symmetric
in the sense that
\[
\langle\partial_s\cdot\chi,\psi\rangle_{\Sigma S^1}=-\langle\chi,\partial_s\cdot\psi\rangle_{\Sigma S^1}
\]
for \(\psi,\chi\in\Gamma(\Sigma S^1)\).
For more details about spinors on \(S^1\) see \cite{MR1682311} and \cite{MR2205370}.\\

For a given curve $\gamma:S^1\to N$, we consider the pull-back bundle $\gamma^{-1}TN$ of $TN$
and twist it with the spinor bundle $\Sigma S^1$.
On this twisted bundle $\Sigma S^1\otimes\gamma^{-1}TN$ there is a metric
induced from the metrics on $\Sigma S^1$ and $\gamma^{-1}TN$.
This induced connection on $\Sigma S^1\otimes\gamma^{-1}TN$
will be denoted by $\tilde{\nabla}$. We will assume that all connections are metric and free of torsion.
Locally, sections of $\Sigma S^1\otimes\gamma^{-1}TN$ can be expressed as
\[
\psi(s)=\psi^i(s)\otimes\frac{\partial}{\partial y^i}(\gamma(s)).
\]
We denote the Dirac operator on $\Sigma S^1$ by $\slashed{\partial}$ and
the Dirac operator on the twisted bundle by $\D$. The untwisted Dirac operator on \(S^1\) is simply given by \(\p=i\nabla\).
In terms of local coordinates \(\D\psi\) can be expressed as
\[
\D\psi=\slashed{\partial}\psi^i\otimes\frac{\partial}{\partial y^i}(\gamma(s))
+\Gamma^i_{jk}(\gamma(s))\nabla\gamma^j\cdot\psi^k(s)\otimes\frac{\partial}{\partial y^i}(\gamma(s)),
\]
where \(\Gamma^i_{jk}\) are the Christoffel symbols on \(N\).
It is easy to see that \(\D\) is self-adjoint with respect to the \(L^2\) norm.
We may now state the energy functional for Dirac-geodesics:
\begin{equation}
E(\gamma,\psi)=\frac{1}{2}\int_{S^1}(|\gamma'|^2+\langle\psi,\D\psi\rangle)ds. 
\end{equation}
In the second term, we use the scalar product on \(\Sigma S^1\otimes\gamma^{-1}TN\).
The Euler-Lagrange equations of the functional \(E(\gamma,\psi)\) are given by (see \cite{MR2262709})
\begin{align}
\label{dirac-harmonic}
\tau(\gamma)&=\cR(\gamma,\psi),\\
\nonumber\D\psi&=0.
\end{align}
Here, \(\tau(\gamma)\) is the tension field of the curve \(\gamma\) and the curvature term on the right hand side 
\(\cR(\gamma,\psi)\) is given by \eqref{curvature-term1}.
Written in coordinates, the Euler-Lagrange equations acquire the form
\begin{align*}
\tau^m(\gamma)-\frac{1}{2}R^m_{~lij}(\gamma)\langle\psi^i,\nabla\gamma^l\cdot\psi^j\rangle_{\Sigma S^1}=&0,\\
\slashed{\partial}\psi^i+\Gamma^i_{jk}(\gamma)\nabla\gamma^j\cdot\psi^k=&0,
\end{align*}
where \(R^m_{~lij}\) denotes the components of the curvature tensor on \(N\).
Solutions of the system (\ref{dirac-harmonic}) are called $\emph{Dirac-geodesics}$ from \(S^1\to N\).
As a limiting case these equations also contain ordinary geodesics, namely for \(\psi=0\).

\section{Regularized Dirac-geodesics and Evolution Equations}
A big problem in the analysis of Dirac-harmonic maps is the fact that the energy functional is unbounded from below.
Moreover, the Dirac operator is only of first order such that tools like the maximum principle are not available.
We try to overcome these difficulties by proposing the following regularization of the energy functional (\(|\tilde{\nabla}\psi|^2=|\D\psi|^2\) since \(M=S^1\))
\begin{equation}
E_\epsilon(\gamma,\psi)=\frac{1}{2}\int_{S^1}(|\gamma'|^2+\langle\psi,\D\psi\rangle+\epsilon|\tilde{\nabla}\psi|^2)ds.
\end{equation}
As a first direct consequence we find
\[
-\frac{1}{8\epsilon}\int_{S^1}|\psi|^2ds\leq E_\epsilon(\gamma,\psi).
\]
Let us now study the critical points of $E_\epsilon(\gamma,\psi)$:
\begin{Prop}
The critical points of \(E_\epsilon(\gamma,\psi)\) are given by
\begin{align}
\label{euler-gamma-regularized}
\tau(\gamma)=&\cR(\gamma,\psi)+\epsilon\cR_c(\gamma,\psi), \\
\label{euler-psi-regularized}
\epsilon\tilde{\Delta}\psi=&\D\psi
\end{align}
and the vector fields \(\cR(\gamma,\psi)\) and \(\cR_c(\gamma,\psi)\) are given by \eqref{curvature-term1} and \eqref{curvature-term2}.
Moreover, \(\tilde{\Delta}\) denotes the connection Laplacian on the bundle \(\Sigma S^1\otimes\gamma^{-1}TN\).
\end{Prop}

\begin{proof}
We start by deriving the Euler-Lagrange equation for the spinor $\psi$. 
Therefore, we consider a variation of $\psi$ with $\gamma$ fixed and $\frac{\tilde{\nabla}\psi_t}{\partial t}\big|_{t=0}=\chi$.
We find
\begin{align*}
\frac{\delta}{\delta\psi}E_\epsilon(\gamma,\psi_t)=&
\frac{1}{2}\int_{S^1}(\langle\chi,\D\psi\rangle
+\langle\psi,\D\chi\rangle+2\epsilon\langle\chi,\tilde{\nabla}^*\tilde{\nabla}\psi\rangle)ds\\
=&\int_{S^1}\langle\chi,\D\psi-\epsilon\tilde{\Delta}\psi\rangle ds.
\end{align*}
To derive the Euler-Lagrange equation for $\gamma$, consider a family of smooth variations of $\gamma$ satisfying 
$\frac{\partial\gamma_t}{\partial t}\big|_{t=0}=\eta$, while keeping the $\psi^i$
in the twisted spinor \(\psi(s)=\psi^i(s)\otimes\frac{\partial}{\partial y^i}(\gamma_t(s))\) fixed.
The variation with respect to $\gamma$ of the following terms has already
been computed in \cite{MR2262709}, p.\ 413:
\begin{align*}
\frac{\delta}{\delta\gamma}\frac{1}{2}\int_{S^1}|\gamma'_t|^2ds=&-\int_{S^1}\langle\tau(\gamma),\eta\rangle ds, \\
\frac{\delta}{\delta\gamma}\frac{1}{2}\int_{S^1}\langle\psi,\D\psi\rangle ds=&\int_{S^1}(\langle\frac{\nabla\psi}{\partial t},\D\psi\rangle
+\langle\cR(\gamma,\psi),\eta\rangle)ds.
\end{align*}
Finally, we compute the variation of the regularizing term, namely
\[
\frac{\partial}{\partial t}\bigg|_{t=0}\frac{\epsilon}{2}\int_{S^1}|\tilde{\nabla}\psi|^2ds=
\epsilon\int_{S^1}\big(\langle\frac{\nabla\psi}{\partial t},\tilde{\nabla}^*\tilde{\nabla}\psi\rangle
+\langle R^E(\partial_t,\partial_s)\psi,\tilde{\nabla}\psi\rangle\big)\big|_{t=0}ds,
\]
where \(R^E\) denotes the curvature tensor on \(E=T^*S^1\otimes\Sigma S^1\otimes\gamma_t^{-1}TN\).
The only curvature contribution arises from the pull-back bundle \(\gamma_t^{-1}TN\) and we compute
\begin{align*}
\langle R^{\gamma_t^{-1}TN}(\partial_t,\partial_s)\psi,\tilde{\nabla}\psi\rangle\big|_{t=0}
&=
\langle R^{N}(\dot{\gamma}_t,\gamma'_t)\psi,\tilde{\nabla}\psi\rangle\big|_{t=0}\\
&=\langle R^{N}(\tilde{\nabla}\psi,\psi)\gamma'_t,\dot{\gamma}_t\rangle\big|_{t=0}\\
&=\langle\cR_c(\gamma,\psi),\eta\rangle.
\end{align*}
Adding up the different contributions, we get
\begin{align*}
\frac{\delta}{\delta\gamma}E_\epsilon(\gamma_t,\psi)=\int_{S^1}(\langle-\tau(\gamma)&+\cR(\gamma,\psi)+\epsilon\cR_c(\gamma,\psi),\eta\rangle \\
&+\langle\D\psi-\epsilon\tilde{\Delta}\psi,\frac{\nabla\psi}{\partial t}\big|_{t=0}\rangle)ds.
\end{align*}
Using the Euler-Lagrange equation for $\psi$, which was deduced before, the result follows.
\end{proof}

Written in local coordinates, the new terms arising from the variation of \(E_\epsilon(\gamma,\psi)\) acquire the following form:
\begin{align*}
\cR_c(\gamma,\psi)=& R^m_{~lij}\frac{\partial}{\partial y^m}\frac{\partial\gamma^l}{\partial s}\langle\nabla^{\Sigma S^1}\psi^i,\psi^j\rangle_{\Sigma S^1}\\
&+ R^m_{~lij}\frac{\partial}{\partial y^m}\Gamma^j_{rt}\frac{\partial\gamma^l}{\partial s}\langle\psi^i,\psi^r\rangle_{\Sigma S^1}\frac{\partial\gamma^t}{\partial s},\\
\tilde{\Delta}\psi=&\Delta^{\Sigma S^1}\psi^i\otimes\frac{\partial}{\partial y^i}
+2\nabla^{\Sigma S^1}\psi^i\otimes\Gamma^k_{ij}\frac{\partial\gamma^j}{\partial s}\frac{\partial}{\partial y^k}\\
&+\psi^i\otimes\Gamma^k_{ij,p}\frac{\partial\gamma^p}{\partial s}\frac{\partial\gamma^j}{\partial s}\frac{\partial}{\partial y^k}
+\psi^i\otimes\Gamma^k_{ij}\frac{\partial^2\gamma^j}{\partial s^2}\frac{\partial}{\partial y^k}\\
&+\psi^i\otimes\Gamma_{ij}^k\Gamma^r_{kt}\frac{\partial\gamma^j}{\partial s}\frac{\partial\gamma^t}{\partial s}\frac{\partial}{\partial y^r}.
\end{align*}
Solutions of the system (\ref{euler-gamma-regularized}), (\ref{euler-psi-regularized}) will be called
\emph{regularized Dirac-geodesics}.

Before turning to the regularized evolution equations, we can give a construction for 
Dirac-geodesics from a closed curve to a Riemannian manifold \(N\), where
we follow the ideas from \cite{MR2262709}, Prop.\ 2.2.
\begin{Prop}
Assume that \(M=S^1\) and \(N\) a compact Riemannian manifold.
If \(\gamma\) is a geodesic and the spinor \(\psi\) is of the form
\begin{equation}
\psi=\partial_s\cdot\chi\otimes\gamma'
\end{equation}
with \(\chi\) being a harmonic spinor, then the pair \((\gamma,\psi)\) 
is a Dirac-geodesic.
\end{Prop}
\begin{proof}
By assumption \(\gamma\) is a geodesic. Hence, we have \(\tau(\gamma)=0\).
By a direct computation, we find that \(\cR(\gamma,\psi)\) is real.
Inserting \(\psi\) as defined above into \(\cR(\gamma,\psi)\) and using that the expression
\[
\overline{\langle\partial_s\cdot\chi,\chi\rangle}_{\Sigma S^1}=\langle\chi, \partial_s\cdot\chi\rangle_{\Sigma S^1}=-\langle \partial_s\cdot\chi,\chi\rangle_{\Sigma S^1}
\]
is purely imaginary, we conclude that \(\cR(\gamma,\psi)\) vanishes.

On the other hand, applying the twisted Dirac operator \(\D\) to the spinor \(\psi\) yields
\begin{align*}
\D\psi=\partial_s\cdot\frac{\tilde{\nabla}}{\partial s}(\partial_s\cdot\chi\otimes\gamma') 
=\partial_s\cdot\p\chi\otimes\gamma'-\chi\otimes\tau(\gamma)
=0
\end{align*}
by assumption, which concludes the proof.
\end{proof}
Let us make two comments about the solution constructed above.
First of all, it is uncoupled, in the sense that \(\tau(\gamma)=0=\cR(\gamma,\psi)\).
Secondly, a harmonic spinor on \(S^1\) only exists for the trivial spin structure \(\sigma_1\) and is given by a constant.

We now turn to studying the \(L^2\)-gradient flow of the regularized functional \(E_\epsilon(\gamma,\psi)\):
\begin{align}
\label{evolution-gamma}
\frac{\partial\gamma_t}{\partial t}=&\tau(\gamma_t)-\cR(\gamma_t,\psi_t)-\epsilon\cR_c(\gamma_t,\psi_t),\\
\label{evolution-psi}
\frac{\tilde{\nabla}\psi_t}{\partial t}=&\epsilon\tilde{\Delta}\psi_t-\D\psi_t.
\end{align}

Our aim is to deform given initial data \((\gamma_0,\psi_0)\) to a regularized Dirac-geodesic by
the gradient flow and then study the limit \(\epsilon\to 0\).
Before we address the general case, let us discuss some illuminating examples.

\theoremstyle{definition}
\begin{Bsp}
Assume that \(M=N=S^1\). Then the evolution equations for \((\gamma_t,\psi_t)\) are given by
\[
\begin{cases}
\partial_t\gamma(s,t)=\partial^2_s\gamma(s,t),\\
\gamma(s,0)=\gamma_0(s),\\
\gamma(0,t)=\gamma(2\pi,t),
\end{cases}
\qquad
\begin{cases}
\partial_t\psi(s,t)=\epsilon\partial^2_s\psi(s,t)-i\partial_s\psi(s,t),\\
\psi(s,0)=\psi_0(s),\\
\psi(0,t)=\pm\psi(2\pi,t).
\end{cases}
\]
The sign in the boundary condition for the spinor depends on the chosen spin structure.
The fundamental solution for the heat equation on $S^1$ can be obtained by a Fourier 
decomposition and is given by
\[
\xi(s,t)=\sum_{k=-\infty}^\infty a_ke^{iks}e^{-k^2t},
\]
with coefficients \(a_k\).

To derive the fundamental solution for the evolution equation for \(\psi(s,t)\),
we make a separation ansatz of the form $\chi(s,t)=A(s)B(t)$ leading to
\[
\frac{\dot B(t)}{B(t)}=C=\frac{-iA'(s)+\epsilon A''(s)}{A(s)}
\]
with a constant \(C\). By \(\lambda_k\) we denote the \(k\)-th eigenvalue of
the Dirac operator on \(S^1\), such that we get
\[
\chi_{\sigma_j}(s,t)=\sum_{k=-\infty}^\infty b_ke^{i\lambda_ks}e^{(\lambda_k-\epsilon\lambda_k^2)t},\qquad j=1,2
\]
with coefficients \(b_k\).
It is known, that for the spin structure \(\sigma_1\), the eigenvalues of the Dirac operator 
are all integer numbers \(\lambda_k=k\),
whereas for the second spin structure \(\sigma_2\) the eigenvalues are given by \(\lambda_k=k+\frac{1}{2}\)
with \(k\in\Z\). Consequently, only $\sigma_1$ admits harmonic spinors.
To incorporate the initial condition at $t=0$, we consider the convolution
\[
\psi(s,t)=\frac{1}{2\pi}\int_0^{2\pi}\psi_0(y)\chi_{\sigma_j}(s-y,t)dy,\qquad j=1,2.
\]
We will analyze the evolution equation for the two spin structures and different choices of $\psi_0(s)$.
To this end, we first fix $\sigma_1$.
\begin{itemize}
\item If $\psi_0(s)=\sum_{k=0}^\infty b_ke^{iks}$, then
\[
\psi(s,t)=\sum_{k=0}^\infty b_ke^{iks}e^{(k-\epsilon k^2)t}.
\]
Without any further restriction on \(b_k,\epsilon\) or \(k\), we cannot make
any statement about \(\psi(s,t)\) as \(t\to\infty\).
\item 
Since \(-\frac{1}{4\epsilon}\leq\lambda_k+\epsilon\lambda_k^2\leq\infty\),
there exists a \(k_0\) such that for all \(k>k_0\) the expression \(\lambda_k+\epsilon\lambda_k^2>0\).
If we now choose the initial spinor $\psi_0(s)=\sum_{k>k_0}^\infty b_ke^{-iks}$, then
\[
\psi(s,t)=\sum_{k>k_0}^\infty b_ke^{-iks}e^{-(k+\epsilon k^2)t}.
\]
In this case the limit $t\to\infty$ exists and we find
\[
\psi_\infty(s)=b_0.
\]
\end{itemize}
As a second step we fix the other spin structure $\sigma_2$.
\begin{itemize}
\item If $\psi_0(s)=\sum_{k=0}^\infty b_ke^{i(k+\frac{1}{2})s}$, 
we again cannot make a general statement as \(t\to\infty\).
\item If $\psi_0(s)=\sum_{k>k_0}^\infty b_ke^{-i(k+\frac{1}{2})s}$ with \(k_0\) as before, then
\[
\psi(s,t)=\sum_{k>k_0}^\infty b_ke^{-i(k+\frac{1}{2})s}e^{-((k+\frac{1}{2})+\epsilon(k+\frac{1}{2})^2)t}.
\]
Again, for this special initial data, the limit \(t\to\infty\) exists,
but the limiting spinor \(\psi_\infty\) will vanish.
\end{itemize}
\end{Bsp}
From the simple example above we learned that 
convergence of the evolution equation for the spinor \(\psi(s,t)\) will depend both on the initial spinor \(\psi_0(s)\)
and the spin structure \(\sigma_j\). Note that the above example could also be studied
for \(\epsilon=0\). The condition on the initial data can be thought of as an APS type condition
as it appears in the context of boundary value problem for Dirac operators \cite{MR0331443}.
Clearly, we realize that if we choose \(\epsilon\geq 1\) we will always get convergence in the evolution equation for \(\psi\).
\\
If we drop the compactness assumption on \(M\), we can find another example
in which the evolution equations for \(\gamma\) and \(\psi\) can be solved explicitly.

\begin{Bsp}
Assume that $M=N=\R$. In this case the evolution equations acquire the form
\[
\begin{cases}
\partial_t\gamma(s,t)=\partial^2_s\gamma(s,t),\\
\gamma(s,0)=\gamma_0(s),
\end{cases}
\qquad
\begin{cases}
\partial_t\psi(s,t)=\epsilon\partial^2_s\psi(s,t)-i\partial_s\psi(s,t),\\
\psi(s,0)=\psi_0(s).
\end{cases}
\]
These equations can be integrated directly. For \(\gamma(s,t)\), we get
the solution to the one-dimensional heat equation
\[
\gamma(s,t)=\frac{1}{\sqrt{4\pi t}} e^{-\frac{s^2}{4t}},
\]
whereas for \(\psi(s,t)\), we find the formal solution
\[
\psi(s,t)=\frac{1}{\sqrt{t}} e^{\frac{1}{4\epsilon}t+\frac{i}{2\epsilon}s-\frac{s^2}{4\epsilon t}}.
\]
Concerning the solution \(\psi(s,t)\), we realize that both limits \(t\to\infty\) and \(\epsilon\to 0\)
are not well defined. This is not surprising since \(M\) is non-compact.
\end{Bsp}

\section{Energy Estimates and Longtime-Existence}
From now on, we study the evolution equations in general.
The existence of a smooth short-time solution up to a time \(T_{max}\) for sufficiently regular initial
data \((\gamma_0,\psi_0)\in C^{2+\alpha}(S^1,N)\times C^{2+\alpha}(S^1,\Sigma S^1\otimes\gamma_0^{-1}TN)\) can be obtained by standard methods 
for semi-linear parabolic partial differential equations. For a detailed proof see Theorem 3.22 in \cite{phd}.
In order to go beyond the shorttime solution we need to derive energy estimates. 
To shorten the notation we introduce the following space
\[
\chi(S^1\times[0,T),N):=C^\infty(S^1\times[0,T),N)\times C^\infty(S^1\times [0,T),\Sigma S^1\otimes\gamma_t^{-1}TN)
\]
and we will use a \(\dot{}\) to represent differentiation with respect to \(t\).
Since we are dealing with a variational problem,
we get bounds in terms of the initial data \((\gamma_0,\psi_0)\).
\begin{Lem}
Let 
\(
(\gamma_t,\psi_t)\in\chi(S^1\times[0,T),N)
\)
be a solution of (\ref{evolution-gamma}) and (\ref{evolution-psi}). 
Then we have for all \(t\in[0,T)\)
\[
E_\epsilon(\gamma_t,\psi_t)+\int_0^T\int_{S^1}\left(\big|\frac{\partial\gamma_t}{\partial t}\big|^2
+\big|\frac{\tilde{\nabla}\psi_t}{\partial t}\big|^2\right)dsdt=E_\epsilon(\gamma_0,\psi_0)
\]
and also
\[
-\frac{1}{8\epsilon}\int_{S^1}|\psi_t|^2ds\leq E_\epsilon(\gamma_t,\psi_t)\leq E_\epsilon(\gamma_0,\psi_0).
\]
\end{Lem}
\begin{proof}
This is a direct consequence of the gradient flow.
\end{proof}

\begin{Bem}
The last Lemma tells us that the regularized energy \(E_\epsilon(\gamma_t,\psi_t)\)
is bounded from below by the \(L^2\)-norm of the spinor \(\psi_t\). Hence, we may expect
that this \(L^2\)-norm will play an important role whenever we will discuss convergence
of the gradient flow.
Moreover, the inequality tells us that the energy is decreasing, which is of course
a general feature of the gradient flow.
\end{Bem}

In the following it turns out to be useful to rescale the \(t\)-parameter in the spinor $\psi(s,t)$ in the following way
\[
\psi(s,t)\to \psi(s,\epsilon t),
\]
such that the pair $(\gamma_t,\psi_t)$ solves the following set of equations
\begin{align}
\label{evolution-gamma-rescaled}
\frac{\partial\gamma_t}{\partial t}=&\tau(\gamma_t)-\cR(\gamma_t,\psi_t)-\epsilon\cR_c(\gamma_t,\psi_t),\\
\label{evolution-psi-rescaled}
\frac{\tilde{\nabla}\psi_t}{\partial t}=&\tilde{\Delta}\psi_t-\frac{1}{\epsilon}\D\psi_t
\end{align}
with initial data \((\gamma(s,0),\psi(s,0))=(\gamma_0(s),\psi_0(s))\).
We will now use the maximum principle to establish pointwise energy estimates.

\begin{Lem}
Let $\psi_t\in C^\infty(S^1\times[0,T),\Sigma S^1\otimes\gamma_t^{-1}TN)$ be a solution of 
(\ref{evolution-psi-rescaled}). Then the norm of the spinor $\psi_t$
satisfies the following estimate:
\begin{equation}
|\psi_t|^2\leq e^{\frac{1}{2\epsilon^2}t}|\psi_0|^2.
\end{equation}
\end{Lem}
\begin{proof}
Due to the rescaling of \(t\) in (\ref{evolution-psi-rescaled}), the norm of $\psi_t$
satisfies the following evolution equation
\[
\frac{\partial}{\partial t}\frac{1}{2}|\psi_t|^2=
\Delta\frac{1}{2}|\psi_t|^2-\frac{1}{\epsilon}\langle\psi_t,\D\psi_t\rangle
-|\tilde{\nabla}\psi_t|^2.
\]
Since $M=S^1$ we have $|\D\psi|^2=|\tilde{\nabla}\psi|^2$.
By Young's inequality we get the estimate
\[
\frac{\partial}{\partial t}\frac{1}{2}|\psi_t|^2\leq\Delta\frac{1}{2}|\psi_t|^2+\frac{1}{4\epsilon^2}|\psi_t|^2.
\]
Now, apply the maximum principle to the function $e^{-\frac{1}{2\epsilon^2}t}|\psi_t|^2$.
\end{proof}

\begin{Cor}
In particular, the last Lemma tells us that if the initial spinor vanishes, that is \(\psi_0=0\), then \(\psi_t=0\) for all \(t\in [0,T)\).
In this case, our set of evolution equations reduces to the heat flow for geodesics studied by Ottarsson \cite{MR834094}.
\end{Cor}

As a next step, we want to derive pointwise estimates 
on the norms of $\gamma'_t$ and $\tilde{\nabla}\psi_t$.

\begin{Prop}
\label{curve-satz-evolution-F}
Let 
\(
(\gamma_t,\psi_t)\in \chi(S^1\times[0,T),N)
\)
be a solution of (\ref{evolution-gamma-rescaled}) and (\ref{evolution-psi-rescaled}).
The function $F_t$ defined by
\begin{equation}
F_t:=\frac{1}{2}(|\gamma'|^2+\epsilon|\tilde{\nabla}\psi_t|^2)
\end{equation}
satisfies the following evolution equation:
\begin{equation}
\label{inequality-F-curve}
\frac{\partial F_t}{\partial t}\leq \Delta F_t+\frac{C}{\epsilon^2}e^{\frac{t}{\epsilon^2}}F_t.
\end{equation}
The constant $C$ depends on $N$ and $\psi_0$.
\end{Prop}
\begin{proof}
Using a standard Bochner technique and remembering the rescaling of \(t\) in the equation for \(\psi_t\),
we find
\begin{align*}
\frac{\partial F_t}{\partial t}
=&\Delta F_t-|\nabla\gamma'_t|^2-\langle\frac{\nabla}{\partial s}\cR(\gamma_t,\psi_t),\gamma'\rangle
-\epsilon\langle\frac{\nabla}{\partial s}\cR_c(\gamma_t,\psi_t),\gamma'\rangle\\
&-\epsilon|\tilde{\nabla}^2\psi_t|^2 +\epsilon\langle\dot\gamma_t,\cR_c(\gamma_t,\psi_t)\rangle-\langle\tilde{\nabla}\D\psi_t,\tilde{\nabla}\psi_t\rangle.
\end{align*}
Note that since $M=S^1$ most of the curvature terms drop out.
A direct computation yields
\[
\langle\cR(\gamma,\psi),\gamma'\rangle_{\gamma^{-1}TN}
=\frac{1}{2}\langle\partial_s\cdot R^N(\gamma',\gamma')\psi,\psi\rangle_{\gamma^{-1}TN}=0.
\]
Therefore, we have
\[
-\langle\frac{\nabla}{\partial s}\cR(\gamma_t,\psi_t),\gamma'\rangle=\langle\cR(\gamma_t,\psi_t),\tau(\gamma_t)\rangle
\]
and similarly for $\cR_c(\gamma_t,\psi_t)$. Applying this identity we find
\begin{align*}
\frac{\partial F_t}{\partial t}
=&\Delta F_t-|\nabla\gamma'_t|^2
-\epsilon^2|\cR_c(\gamma_t,\psi_t)|^2-\epsilon|\tilde{\nabla}^2\psi_t|^2-\langle\tilde{\nabla}\D\psi_t,\tilde{\nabla}\psi_t\rangle\\
&+\langle 2\epsilon\cR_c(\gamma_t,\psi_t)+\cR(\gamma_t,\psi_t),\tau(\gamma_t)\rangle
-\epsilon\langle\cR_c(\gamma_t,\psi_t),\cR(\gamma_t,\psi_t)\rangle.
\end{align*}
As a next step, we use the estimate
\begin{align*}
-&|\nabla\gamma'_t|^2-\epsilon^2|\cR(\gamma_t,\psi_t)|^2-\epsilon\langle\cR_c(\gamma_t,\psi_t),\cR(\gamma_t,\psi_t)\rangle\\
&\hspace{0.4cm}+\langle 2\epsilon\cR_c(\gamma_t,\psi_t)+\cR(\gamma_t,\psi_t),\tau(\gamma_t)\rangle \\
&\leq
\frac{1}{4}|\cR(\gamma_t,\psi_t)+2\epsilon\cR_c(\gamma_t,\psi_t)|^2
-\epsilon^2|\cR_c(\gamma_t,\psi_t)|^2-\epsilon\langle\cR_c(\gamma_t,\psi_t),\cR(\gamma_t,\psi_t)\rangle\\
&=\frac{1}{4}|\cR(\gamma_t,\psi_t)|^2
\end{align*}
and apply Young's inequality again
\[
-\langle\tilde\nabla\D\psi_t,\tilde\nabla\psi_t\rangle-\epsilon|\tilde\nabla^2\psi_t|^2\leq
\frac{1}{4\epsilon}|\tilde\nabla\psi_t|^2.
\]
Finally, we calculate
\begin{align*}
\frac{\partial F_t}{\partial t}\leq & \Delta F_t+\frac{1}{4}|\cR(\gamma_t,\psi_t)|^2+\frac{1}{4\epsilon}|\tilde\nabla\psi_t|^2\\
\leq&  \Delta F_t+C\big(|\psi_t|^4|\gamma'_t|^2+\frac{1}{\epsilon}|\tilde\nabla\psi_t|^2\big)\\
\leq&  \Delta F_t+C\big(e^{\frac{1}{\epsilon^2}t}|\gamma'_t|^2+\frac{1}{\epsilon}|\tilde\nabla\psi_t|^2\big)\\
\leq&  \Delta F_t+Ce^{\frac{1}{\epsilon^2}t}\big(|\gamma'_t|^2+\frac{1}{\epsilon}|\tilde\nabla\psi_t|^2\big)\\
\leq&  \Delta F_t+Ce^{\frac{1}{\epsilon^2}t}\frac{1}{\epsilon^2}\big(|\gamma'_t|^2+\epsilon|\tilde\nabla\psi_t|^2\big)\\
\leq&  \Delta F_t+\frac{C}{\epsilon^2}e^{\frac{1}{\epsilon^2}t}F_t,
\end{align*}
where we used the fact that $\epsilon<1$.
\end{proof}

\begin{Cor}
Let 
\(
(\gamma_t,\psi_t)\in\chi(S^1\times[0,T),N)
\)
be a solution of (\ref{evolution-gamma-rescaled}) and (\ref{evolution-psi-rescaled}).
Then we have for all \(t\in [0,T)\)
\begin{equation}
F_t=\frac{1}{2}(|\gamma'_t|^2+\epsilon|\tilde{\nabla}\psi_t|^2)\leq C e^{Ce^\frac{t}{\epsilon^2}}.
\end{equation}
\end{Cor}
The constant $C$ depends on \(N,\epsilon,\psi_0,\gamma'_0\) and \(\tilde{\nabla}\psi_0\).
\begin{proof}
The evolution equation for \(F_t\) (\ref{inequality-F-curve}) can be rewritten as
\[
\frac{\partial}{\partial t}(e^{-Ce^{\frac{t}{\epsilon^2}}}F_t)\leq \Delta(e^{-Ce^{\frac{t}{\epsilon^2}}}F_t).
\]
Applying the maximum principle, we find
\[
F_t\leq e^{-C}e^{Ce^{\frac{t}{\epsilon^2}}}F_{0}.
\]
Rearranging the constants yields the result.
\end{proof}

Having obtained pointwise estimates for $|\gamma'_t|^2$ and $|\tilde{\nabla}\psi_t|^2$,
we can now derive estimates on the $t$-derivatives of $\gamma_t$ and $\psi_t$.
\begin{Prop}
\label{curve-satz-evolution-G}
Let the pair
\(
(\gamma_t,\psi_t)\in\chi(S^1\times[0,T),N)
\)
be a solution of (\ref{evolution-gamma-rescaled}) and (\ref{evolution-psi-rescaled}).
The quantity defined by
\begin{equation}
G_t:=\frac{1}{2}\big(|\frac{\partial\gamma_t}{\partial t}|^2+|\frac{\tilde{\nabla}\psi_t}{\partial t}|^2\big)
\end{equation}
satisfies
\begin{equation}
\frac{\partial G_t}{\partial t}\leq \Delta G_t+Z(t)G_t
\end{equation}
with the function $Z(t)$ depending on $\epsilon,|\psi_t|,|\gamma_t|$, and $|\tilde{\nabla}\psi_t|$.
\end{Prop}

\begin{proof}
By a direct but lengthy calculation, taking into account the rescaling of \(t\) in the evolution equation for \(\psi_t\) and the fact that $M=S^1$, we find
\begin{align*}
\frac{\partial G_t}{\partial t} =&\Delta G_t 
-|\nabla\dot{\gamma}_t|^2-|\tilde{\nabla}\frac{\tilde{\nabla}\psi_t}{\partial t}|^2
+\langle R^N(\gamma'_t,\dot{\gamma}_t)\gamma'_t,\dot{\gamma}_t\rangle
-\langle\frac{\nabla}{\partial t} \cR(\gamma_t,\psi_t),\dot{\gamma}_t\rangle \\
&-\epsilon\langle\frac{\nabla}{\partial t}\cR_c(\gamma_t,\psi_t),\dot{\gamma}_t\rangle
-\frac{1}{\epsilon}\langle\frac{\tilde{\nabla}\psi_t}{\partial t},\D\frac{\tilde{\nabla}\psi_t}{\partial t}\rangle-\frac{1}{\epsilon}\langle\partial_s\cdot R^N(\dot{\gamma}_t,\gamma'_t)\psi_t,\frac{\tilde{\nabla}\psi_t}{\partial t}\rangle\\
&+\langle R^N(\dot{\gamma}_t,\gamma'_t)\tilde{\nabla}\psi_t,\frac{\tilde{\nabla}\psi_t}{\partial t}\rangle
+\langle \tilde{\nabla}(R^N(\dot{\gamma}_t,\gamma'_t)\psi_t),\frac{\tilde{\nabla}\psi_t}{\partial t}\rangle.
\end{align*}
Again, we have to estimate all terms on the right hand side.
First of all, we estimate the term containing the curvature of $\gamma^{-1}TN$ as
\[
\langle R^N(\gamma'_t,\dot{\gamma}_t)\gamma'_t,\dot{\gamma}_t\rangle
\leq C|\dot{\gamma}_t|^2|\gamma'_t|^2:=I_{1}.
\]
For the next term we calculate  
\begin{align*}
\frac{\nabla}{\partial t}\cR(\gamma_t,\psi_t)=&
\frac{1}{2}(\nabla_{\dot{\gamma}_t}R^N)(\partial_s\cdot\psi_t,\psi_t)\gamma'_t
+R^N(\partial_s\cdot\psi_t,\frac{\tilde{\nabla}\psi_t}{\partial t})\gamma'_t\\
&+\frac{1}{2}R^N(\partial_s\cdot\psi_t,\psi_t)\frac{\nabla}{\partial t}\gamma'_t
\end{align*}
and make the following estimate
\begin{align*}
|\langle\frac{\nabla}{\partial t}\cR(\gamma_t,\psi_t),\dot{\gamma}_t\rangle|
&\leq C\big(|\dot{\gamma}_t|^2|\gamma'_t||\psi_t|^2+
|\dot{\gamma}_t||\frac{\tilde{\nabla}\psi_t}{\partial t}||\psi_t||\gamma'_t|
+|\dot{\gamma}_t||\psi_t|^2|\frac{\nabla}{\partial t}\gamma'_t|\big)\\
&:=I_2+I_3+I_4.
\end{align*}
To take care of the next term, we first compute
\begin{align*}
\frac{\nabla}{\partial t}\cR_c(\gamma_t,\psi_t)=&
(\nabla_{\dot{\gamma}_t}R^N)(\tilde{\nabla}\psi_t,\psi_t)\gamma'_t
+R^N(\frac{\tilde{\nabla}}{\partial t}\tilde{\nabla}\psi_t,\psi_t)\gamma'_t\\
&+R^N(\tilde{\nabla}\psi_t,\frac{\tilde{\nabla}\psi_t}{\partial t})\gamma'_t
+R^N(\psi_t,\tilde{\nabla}\psi_t)\frac{\nabla}{\partial t}\gamma'_t.
\end{align*}
Changing covariant derivatives once more
\[
\frac{\tilde{\nabla}}{\partial t}\tilde{\nabla}\psi_t=\tilde{\nabla}\frac{\tilde{\nabla}}{\partial t}\psi_t+R^N(\dot{\gamma}_t,\gamma'_t)\psi_t
\]
we may estimate
\begin{align*}
|\langle\frac{\nabla}{\partial t}\cR_c(\gamma_t,\psi_t),\dot{\gamma}_t\rangle|
\leq C\big(&|\dot{\gamma}_t|^2|\psi_t||\gamma'_t||\tilde{\nabla}\psi_t|
 +|\dot{\gamma}_t||\psi_t||\tilde{\nabla}\frac{\tilde{\nabla}\psi}{\partial t}||\gamma'_t| 
 +|\dot{\gamma}_t|^2|\psi_t||\gamma'_t|^2 \\
& +|\dot{\gamma}_t||\frac{\tilde{\nabla}\psi_t}{\partial t}||\tilde{\nabla}\psi_t||\gamma'_t|
+|\dot{\gamma}_t||\psi_t||\tilde{\nabla}\psi_t||\frac{\nabla}{\partial t}\gamma'_t|\big) \\
:=I_5&+I_6+I_7+I_8+I_9.
\end{align*}
Having estimated the terms from the evolution equation for $\gamma_t$, we now deal 
with the contributions originating from the evolution equation for \(\psi_t\)
\begin{align*}
\frac{1}{\epsilon}\big|\langle\frac{\tilde{\nabla}\psi_t}{\partial t},\D\frac{\tilde{\nabla}\psi_t}{\partial t}\rangle\big|
&\leq \frac{2}{\epsilon^2}|\frac{\tilde{\nabla}\psi_t}{\partial t}|^2
+\frac{1}{8}|\tilde{\nabla}\frac{\tilde{\nabla}\psi_t}{\partial t}|^2:=I_{10}+I_{11},\\
\big|\langle\partial_s\cdot R^N(\dot{\gamma}_t,\gamma'_t)\psi_t,\frac{\tilde{\nabla}\psi_t}{\partial t}\rangle\big|
&\leq C|\dot{\gamma}_t||\frac{\tilde{\nabla}\psi_t}{\partial t}||\psi_t||\gamma'_t|:=I_{12},\\
|\langle R^N(\dot{\gamma}_t,\gamma'_t)\tilde{\nabla}\psi_t,\frac{\tilde{\nabla}\psi_t}{\partial t}\rangle|
&\leq C|\dot{\gamma}_t||\frac{\tilde{\nabla}\psi_t}{\partial t}||\tilde{\nabla}\psi_t||\gamma'_t|:=I_{13}.
\end{align*}
To estimate the last term in the evolution equation for \(G_t\) we compute
\begin{align*}
\frac{\tilde{\nabla}}{\partial s}(R^N(\dot{\gamma}_t,\gamma'_t)\psi_t)
=& (\nabla_{\gamma'_t}R^N)(\dot{\gamma}_t,\gamma'_t)\psi_t+R^N(\frac{\nabla}{\partial s}\dot{\gamma}_t,\gamma'_t)\psi_t \\
&+R^N(\dot{\gamma}_t,\tau(\gamma_t))\psi_t +R^N(\dot{\gamma}_t,\gamma'_t)\tilde{\nabla}\psi_t.
\end{align*}
The first, second and fourth term can easily be estimated as
\begin{align*}
|\langle(\nabla_{\gamma_t'}R^N)(\dot{\gamma}_t,\gamma'_t)\psi_t,\frac{\tilde{\nabla}\psi_t}{\partial t}\rangle| &\leq  C|\frac{\tilde{\nabla}\psi_t}{\partial t}||\dot{\gamma}_t||\gamma'_t|^2|\psi_t|, \\
|\langle R^N(\frac{\nabla}{\partial s}\dot{\gamma}_t,\gamma'_t)\psi_t,\frac{\tilde{\nabla}\psi_t}{\partial t}\rangle|&\leq C|\frac{\tilde{\nabla}\psi_t}{\partial t}||\nabla\dot{\gamma}_t||\gamma'_t||\psi_t|,\\
|\langle R^N(\dot{\gamma}_t,\gamma'_t)\tilde{\nabla}\psi_t,\frac{\tilde{\nabla}\psi_t}{\partial t}\rangle|&\leq C |\dot{\gamma}_t||\frac{\tilde{\nabla}\psi_t}{\partial t}||\gamma'_t||\tilde{\nabla}\psi_t|.
\end{align*}
To take care of the third term, we use the evolution equation for $\gamma_t$ and obtain
\[
R^N(\dot{\gamma}_t,\tau(\gamma_t))\psi_t
=R^N(\dot{\gamma}_t,\cR(\gamma_t,\psi_t))\psi_t+\epsilon R^N(\gamma_t,\cR_c(\gamma_t,\psi_t))\psi_t.
\]
This allows us to derive the following estimate
\begin{align*}
|\langle R^N(\dot{\gamma}_t&,\tau(\gamma_t))\psi_t,\frac{\tilde{\nabla}\psi_t}{\partial t}\rangle|\\
&\leq |\langle R^N(\dot{\gamma}_t,\cR(\gamma_t,\psi_t))\psi_t,\frac{\tilde{\nabla}\psi_t}{\partial t}\rangle|
+\epsilon|\langle R^N(\dot{\gamma}_t,\cR_c(\gamma_t,\psi_t))\psi_t,\frac{\tilde{\nabla}\psi_t}{\partial t}\rangle|\\
&\leq C\big(|\frac{\tilde{\nabla}\psi_t}{\partial t}||\dot{\gamma}_t||\gamma'_t||\psi_t|^3
+\epsilon|\dot{\gamma}_t||\tilde{\nabla}\psi_t||\gamma'_t||\psi_t|^2|\frac{\tilde{\nabla}\psi_t}{\partial t}|\big)
\end{align*}
and, finally, we have
\begin{align*}
|\langle\tilde{\nabla}(R^N(\dot{\gamma}_t,&\gamma'_t)\psi_t),\frac{\tilde{\nabla}\psi_t}{\partial t})\rangle|\\
&\leq C\big(|\dot{\gamma}_t||\frac{\tilde{\nabla}\psi_t}{\partial t}||\gamma'_t|^2|\psi_t|
+|\nabla\dot{\gamma}_t||\frac{\tilde{\nabla}\psi_t}{\partial t}||\gamma'_t||\psi_t|
+|\dot{\gamma}_t||\frac{\tilde{\nabla}\psi_t}{\partial t}||\gamma'_t||\psi_t|^3\\
&\hspace{1cm}+\epsilon|\dot{\gamma}_t||\frac{\tilde{\nabla}\psi_t}{\partial t}||\gamma'_t||\psi_t|^2|\tilde{\nabla}\psi_t|
+|\dot{\gamma}_t||\frac{\tilde{\nabla}\psi_t}{\partial t}||\gamma'_t||\tilde{\nabla}\psi_t|\big)\\
&:=I_{14}+I_{15}+I_{16}+I_{17}+I_{18}.
\end{align*}
Collecting all the estimates, the function \(G_t\) satisfies
\begin{equation}
\label{G-inequality-a}
\frac{\partial G_t}{\partial t}\leq \Delta G_t
-|\nabla\dot{\gamma}_t|^2-|\tilde\nabla\frac{\tilde{\nabla}\psi_t}{\partial t}|^2
+\sum_{j=1}^{18}I_j.
\end{equation}
We want to use $-|\nabla\dot{\gamma}_t|^2-|\tilde{\nabla}\frac{\tilde{\nabla}\psi_t}{\partial t}|^2$ in order 
to control $I_4,I_6,I_9,I_{11}$ and $I_{15}$, more precisely
\begin{align*}
&-|\nabla\dot{\gamma}_t|^2-|\nabla\frac{\tilde{\nabla}\psi_t}{\partial t}|^2
+I_4+I_6+I_9+I_{11}+I_{15}\\
&\leq C\big(|\psi_t|^4|\dot{\gamma}_t|^2+\epsilon^2|\psi|^2|\dot{\gamma}_t|^2|\gamma'_t|^2
+\epsilon^2|\dot{\gamma}_t|^2|\tilde{\nabla}\psi_t|^2|\psi_t|^2
+|\gamma'_t|^2|\psi_t|^2|\frac{\tilde{\nabla}\psi_t}{\partial t}|^2\big).
\end{align*}
Hence, we can write
\begin{equation}
\label{G-inequality-b}
\sum_{j=1}^{18}I_j-|\nabla\dot{\gamma}_t|^2-|\tilde{\nabla}\frac{\tilde{\nabla}\psi_t}{\partial t}|^2 \leq
A_t|\dot{\gamma}_t|^2 
+B_t|\frac{\tilde{\nabla}\psi_t}{\partial t}|^2 
+C_t|\dot{\gamma}_t||\frac{\tilde{\nabla}\psi_t}{\partial t}|
\end{equation}
with the terms
\begin{align*}
A_t &=C\big(|\psi_t|^4+\epsilon^2|\psi_t|^2|\gamma'_t|^2+\epsilon^2|\psi_t|^2|\tilde{\nabla}\psi_t|^2+|\gamma'_t|^2+|\gamma'_t||\psi_t|^2 \\
&\hspace{1cm}+\epsilon|\psi_t||\gamma'_t||\tilde{\nabla}\psi_t|\big),\\
B_t &=C\big(|\psi_t|^2|\gamma'_t|^2+\frac{1}{\epsilon^2}\big),\\
C_t &=C\big(|\psi_t||\gamma'_t|+\epsilon|\tilde{\nabla}\psi_t||\gamma'_t|+\frac{1}{\epsilon}|\psi_t||\gamma'_t|+|\tilde{\nabla}\psi_t||\gamma'_t|+|\gamma'_t|^2|\psi_t|+|\gamma'_t||\psi_t|^3\\
&\hspace{1cm}+|\gamma'_t||\psi_t|^2|\tilde{\nabla}\psi_t|\big).
\end{align*}
By the bounds on \(\psi_t,\gamma'_t\), and \(\tilde{\nabla}\psi_t\) all terms appearing in $A_t,B_t$ and $C_t$ can be controlled.
We combine them into one function \(Z(t)\).
Together with \eqref{G-inequality-a} and \eqref{G-inequality-b} this completes the proof.
\end{proof}

\begin{Cor}
Let the pair
\(
(\gamma_t,\psi_t)\in\chi(S^1\times[0,T),N)
\)
be a solution of (\ref{evolution-gamma-rescaled}) and (\ref{evolution-psi-rescaled}).
For all $(s,t)\in S^1\times [0,T)$, we have the following estimate
\begin{equation}
G_t=\frac{1}{2}\big(|\frac{\partial\gamma_t}{\partial t}|^2+|\frac{\tilde{\nabla}\psi_t}{\partial t}|^2\big)\leq C f(t),
\end{equation}
where \(f(t)\) only depends on \(t\) and is finite for finite values of \(t\).
The constant \(C\) depends on $N,\epsilon,\gamma'_0,\psi_0\) and \(\tilde{\nabla}\psi_0$.
\end{Cor}
\begin{proof}
We use the inequality for $G_t$ and apply the maximum principle
to $e^{-\int_0^TZ(\tau)d\tau}G_t$.
The function \(Z(t)\) can be explicitly expressed in terms of exponentials
and double exponentials. These are bounded for finite \(t\) and, consequently,
the integral \(\int_0^TZ(\tau)d\tau\) is also finite.
\end{proof}

Using the estimates just derived, we now turn to establishing a stability and uniqueness result
that will be needed in order to obtain the existence of a unique long-time solution of the evolution equations.
To this end, it is necessary to get an estimate on expressions like \(\gamma_1-\gamma_2\) and 
\(\psi_1-\psi_2\). We therefore apply the embedding theorem of Nash to isometrically embed the manifold \(N\)
into some \(\R^q\) via the map \(\iota\colon N\to\R^q\).
Then, \(u=\iota\circ\gamma\colon S^1\to \R^q\) can be thought of as a vector-valued function.
From the point of view of the manifold \(N\), the spinor \(\psi\) behaves like a tangent vector. 
This is the reason why we have to use the differential of the map \(\iota\) if we want to define
a spinor along the composite curve \(u=\iota\circ\gamma\). More precisely, if $\psi$ is a spinor along $\gamma$, 
then $\hat\psi$ is a spinor along $u$. Both are related by
\begin{equation}
\hat\psi=d\iota(\psi).
\end{equation}
The vector spinor \(\psi\) then turns into a vector of usual spinors \(\hat\psi=(\psi^1,\ldots,\psi^q)\)
with \(\psi^i\in\Gamma(\Sigma S^1),~i=1,\ldots,q\). The condition that \(\psi\) is along the curve \(\gamma\) is encoded by
\[
\sum_{i=1}^q\nu_i\psi^i=0 \qquad \textrm{for any normal vector } \nu_i\in\R^q \textrm{ at } \gamma(s).
\]
Regarding the evolution equations (\ref{evolution-gamma-rescaled}) and (\ref{evolution-psi-rescaled}) we are now in the situation that
\(u\colon S^1\times[0,T)\to\R^q\) and \(\hat\psi\colon S^1\times[0,T)\to\Sigma S^1\otimes T\R^q\).
More precisely, \(u\) satisfies the following equation:
\begin{align}
\label{evolution-u-Rq}
\nonumber
\left(\frac{\partial}{\partial t}-\Delta\right)u=&-\sff_u(u',u')
-P(\sff(u',\partial_s\cdot\hat\psi),\hat\psi)
+\epsilon P(\sff(u',\hat\psi),\nabla\hat\psi)\\
 &-\epsilon P(\sff(u',\nabla\hat\psi),\hat\psi)
-\epsilon B(u',\psi,u',\hat\psi)
\end{align}
with the initial condition \(u_0=\iota(\gamma_0)\). 
For the spinor \(\hat\psi\in\Gamma(\Sigma S^1\otimes T\R^q)\), we get
\begin{align}
\label{evolution-psi-Rq}
\nonumber
\bigg(\frac{\nabla}{\partial t}-\Delta\bigg)\hat\psi=&-\frac{1}{\epsilon}\p\hat\psi
+\frac{1}{\epsilon}\sff(u',\partial_s\cdot\hat\psi)+\sff(\dot{u},\hat\psi)
-2\sff(u',\nabla\hat\psi)\\
 &-(\nabla\sff)(u',\hat\psi))-\sff(\tau(u),\hat\psi)
\end{align}
with the initial condition \(\hat\psi_0=d\iota(\psi_0)\).
Here, \(\sff\) denotes the second fundamental form of the embedding, \(P\) is the shape operator and the quantity \(B\) is given by
\begin{align*}
B_u(u',&\psi,u',\psi)=\\
&\Gamma^m_{kl}(P(\sff_u(\partial_{y^i},\partial_{y^m}),\partial_{y^j})
-P(\sff_u(\partial_{y^i},\partial_{y^j}),\partial_{y^m}))
\frac{\partial u^i}{\partial s}\frac{\partial u^k}{\partial s}\psi^{l}\psi^{j}.
\end{align*}
For more details, see Chapter 3 in \cite{phd}.	
To shorten the notation, we will omit the superscript \(\hat{}\) at the spinors from now on.	

After these preparations we are now able to give the following
\begin{Prop}[Stability and uniqueness]
\label{theorem-curve-uniqueness}
Let \((\gamma_t,\psi_t)\) and \((\tilde{\gamma}_t,\tilde{\psi}_t)\) be
smooth solutions of (\ref{evolution-gamma-rescaled}) and (\ref{evolution-psi-rescaled}).
The spinor \(\psi_t\in\Gamma(\Sigma S^1\otimes\gamma_t^{-1}TN)\) is defined along the curve \(\gamma_t\) and the spinor \(\tilde{\psi}_t\in\Gamma(\Sigma S^1\otimes\tilde{\gamma}_t^{-1}TN)\) along the curve \(\tilde{\gamma_t}\).
If the initial data coincides, i.e. \((\gamma_0,\psi_0)=(\tilde{\gamma}_0,\tilde{\psi}_0)\), 
then we have $(\gamma_t,\psi_t)=(\tilde{\gamma}_t,\tilde{\psi}_t)$ throughout $S^1\times[0,T)$.
\end{Prop}

\begin{proof}
We apply the Nash embedding theorem and set \(u:=\iota\circ\gamma, v:=\iota\circ\tilde{\gamma}\) and 
\(\hat\psi=d\iota(\psi), \xi=d\iota(\tilde{\psi	})\). Again, we omit the superscript \(\hat{}\) at the spinor \(\psi\).
We regard both \(u,v\) as vector valued functions in \(\R^q\), 
i.e. \(u,v\colon S^1\times [0,T)\to\iota(N)\subset\R^q\), and
the spinors as \(\psi,\xi\colon S^1\times [0,T)\to\Sigma S^1\otimes T\R^q\). 
First of all, we define functions \(h_1,h_2\) 
\[
h_1\colon S^1\times[0,T)\to\R^q,\qquad h_2\colon S^1\times[0,T)\to\Sigma S^1\otimes T\R^q
\]
by
\[
h_1=u-v,\qquad h_2=\psi-\xi.
\]
Using the evolution equation for \(u\) derived in (\ref{evolution-u-Rq}), we calculate
\begin{align*}
\frac{\partial}{\partial t}\frac{1}{2}|h_1|^2=&\Delta\frac{1}{2}|h_1|^2-|h'_1|^2
+\langle h_1, \sff_u(u',u')-\sff_v(v',v')\rangle\\
&+\langle h_1,P(\sff_u(u',\partial_s\cdot\psi),\psi)-P(\sff_v(v',\partial_s\cdot\xi),\xi)\rangle\\
&+\epsilon\langle h_1,P(\sff_u(u',\nabla\psi),\psi)-P(\sff_v(v',\nabla\xi),\xi)\rangle\\
&+\epsilon\langle h_1,P(\sff_u(u',\psi),\nabla\psi)-P(\sff_v(v',\xi),\nabla\xi)\rangle\\
&+\epsilon\langle h_1,B_u(u',\psi,u',\psi)-B_v(v',\xi,v',\xi)\rangle.
\end{align*}
We want to estimate the right hand side in terms of the functions \(h_1\) and \(h_2\).
To this end, we use the bounds on \(\gamma'_t,\tilde{\nabla}\psi_t\), and \(\psi_t\)
derived before. Rearranging the second fundamental forms, 
\[
\sff_u(u',u')-\sff_v(v',v')=
(\sff_u-\sff_v)(u',u')+\sff_v(u'-v',u')+\sff_v(v',u'-v'),
\]
and applying the mean value theorem, we find
\[
|\langle \sff_u(u',u')-\sff_v(v',v'),h_1\rangle|\leq C(|h_1|^2+|h_1'||h_1|).
\]
We rewrite
\begin{align*}
P(\sff_u(u'&,\partial_s\cdot\psi),\psi)-P(\sff_v(v',\partial_s\cdot\xi),\xi)
\\&=P(\sff_{u-v}(u',\partial_s\cdot\psi),\psi)
+P(\sff_v(u'-v',\partial_s\cdot\psi),\psi)\\
&\hspace{0.4cm}+P(\sff_v(v',\partial_s\cdot(\psi-\xi)),\psi)
+P(\sff_v(v',\partial_s\cdot\psi),\psi-\xi)
\end{align*}
and estimate again
\begin{align*}
|\langle h_1,P(\sff_u(u',\partial_s&\cdot\psi),\psi)-P(\sff_v(v',\partial_s\cdot\xi),\xi)\rangle| \\
&\leq C(|u'||h_1|^2|\psi|^2+|h_1'||\psi|^2|h_1|+|v'||\psi||h_2||h_1|)\\
&\leq C(||h_1|^2+|h_1'||h_1|+|h_2||h_1|).
\end{align*}
Again, we rewrite
\begin{align*}
P(\sff_u(u'&,\psi),\nabla\psi)-P(\sff_v(v',\xi),\nabla	\xi)
\\&=P(\sff_{u-v}(u',\psi),\nabla\psi)
+P(\sff_v(u'-v',\psi),\nabla\psi)\\
&\hspace{0.4cm}+P(\sff_v(v',\psi-\xi),\nabla\psi)
+P(\sff_v(v',\psi),\nabla(\psi-\xi))
\end{align*}
and estimate 
\begin{align*}
&|\langle h_1,P(\sff_u(u',\psi),\nabla\psi)-P(\sff_v(v',\xi),\nabla\xi)\rangle| \\
&\leq C(|u'||h_1|^2|\psi||\nabla\psi|+|h_1'||\psi||\nabla\psi||h_1|+|v'||\nabla\psi||h_2||h_1|+|v'||\psi||\nabla h_2||h_1|)\\
&\leq C(|h_1|^2+|h_1'||h_1|+|h_2||h_1|+|\nabla h_2||h_1|).
\end{align*}
The term
\[
\langle h_1,P(\sff_u(u',\nabla\psi),\psi)-P(\sff_v(v',\nabla\xi),\xi)\rangle\\
\]
can be treated by the same methods and estimated like the previous one.
Finally, we rewrite 
\begin{align*}
B_u(u',\psi&,u',\psi)-B_v(v',\xi,v',\xi)\\
&=B_{u-v}(u',\psi,u',\psi)+B_v(u'-v',\psi,u',\psi)+B_v(v',\psi-\xi,u',\psi)\\
&\hspace{0.5cm}+B_v(v',\xi,u'-v',\psi)+B_v(v',\xi,v',\psi-\xi)
\end{align*}
such that we can estimate 
\begin{align*}
|\langle h_1,B_u(u'&,\psi,u',\psi)-B_v(v',\xi,v',\xi)\rangle|\\
&\leq C(|h_1|^2|\psi|^2|u'|^2+|h_1||h_1'||\psi|^2|u'|+|u'||\psi||v'||h_2||h_1|\\
&\hspace{0.6cm}+|h_1||v'|\xi||\psi||h_1'|+|h_1||v'|^2|\xi||h_2|)\\
&\leq C(|h_1|^2+|h_1||h_1'|+|h_2||h_1|).
\end{align*}

Collecting all the terms and applying Young's inequality, we find that the norm of \(h_1\) satisfies 
\[
\frac{\partial}{\partial t}\frac{1}{2}|h_1|^2\leq \Delta\frac{1}{2}|h_1|^2 -\frac{1}{2}|h'_1|^2+\frac{1}{2}|\nabla h_2|^2+
C(|h_1|^2+|h_2|^2).
\]
Now, we turn to the function \(h_2\). Using the evolution equation (\ref{evolution-psi-Rq}) a direct computation yields 
\begin{align*}
\frac{\partial}{\partial t}\frac{1}{2}|h_2|^2=\Delta \frac{1}{2}|h_2|^2&-|\nabla h_2|^2
-\frac{1}{\epsilon}\langle\p h_2,h_2\rangle \\
&+\langle h_2,(\nabla\sff_u)(u',\psi)-(\nabla\sff_v)(v',\xi)\rangle. 
\end{align*}
The other terms involving the second fundamental form \(\sff\) vanish since we have \(\sff\perp\psi\). 
The last term in the equation for \(h_2\) can be rewritten as
\begin{align*}
(\nabla\sff_u)(u',\psi)-(\nabla\sff_v)&(v',\xi)=
(\nabla\sff_{h_1})(u',\psi)
+(\nabla\sff_v)(h_1',\psi)
+(\nabla\sff_v)(v',h_2)
\end{align*}
and we may estimate 
\begin{align*}
|\langle h_2,(\nabla\sff_u)(u',\psi)&-(\nabla\sff_v)(v',\xi)\rangle|\\
&\leq C(|h_2||h_1'||u'||\psi|+|v'||h_1'||\psi||h_2|+|v'|^2|h_2|^2)\\
&\leq C(|h_2|^2+|h_1'||h_2|+|h_2|^2).
\end{align*}
Hence, after applying Young's inequality the norm of \(h_2\) satisfies
\[
\frac{\partial}{\partial t}\frac{1}{2}|h_2|^2\leq\Delta\frac{1}{2}|h_2|^2-\frac{1}{2}|\nabla h_2|^2+\frac{1}{2}|h'_1|^2+C(|h_1|^2+|h_2|^2).
\]
Finally, we define the function \(h\colon S^1\times[0,T)\to\R\) as \(h:=\frac{1}{2}(|h_1|^2+|h_2|^2)\).
One can think of \(h\) as the norm of \((h_1,h_2)\).
Clearly, \(h\) satisfies the following inequality
\[
\frac{\partial h}{\partial t}\leq \Delta h+Ch.
\]
By the maximum principle we get 
\[
\max_{S^1\times[0,T)} h(s,t)\leq\max_{S^1} h(x,0)e^{Ct},
\]
but by assumption \(h(s,0)=0\). Thus, we have \(u=v\) and also \(\psi=\xi\)
throughout \(S^1\times[0,T)\).
\end{proof}

In the next Proposition, we improve the regularity of the pair \((\gamma_t,\psi_t)\)
by application of the classical Schauder estimates, see for example \cite{MR2309679}, Chapters 6 and 7. 
In the one-dimensional case considered here, elliptic Schauder theory is of course not more than just integrating the right hand side.

\begin{Prop}[Schauder theory]
\label{curve-prop-schauder}
Let \((\gamma_t,\psi_t)\in\chi(S^1\times[0,T),N)\) be a solution of (\ref{evolution-gamma-rescaled}) and (\ref{evolution-psi-rescaled}).
Then for any $0<\alpha<1$, there exists a positive number $C$ such that
\begin{align}
|\gamma(\cdot,t)|_{C^{2+\alpha}(S^1,N)}+\big|\frac{\partial\gamma}{\partial t}(\cdot,t)\big|_{C^\alpha(S^1,N)}&\leq C,\\
|\psi(\cdot,t)|_{C^{2+\alpha}(S^1,\Sigma S^1\otimes\gamma_t^{-1}TN)}+\big|\frac{\tilde{\nabla}\psi}{\partial t}(\cdot,t)\big|_{C^{\alpha}(S^1,\Sigma S^1\otimes\gamma_t^{-1}TN)}&\leq C
\end{align}
hold for all $t\in [0,T)$, where both constants depend on \(N,\epsilon,\alpha,T,\psi_t,\gamma'_t\) and \(\tilde{\nabla}\psi_t\).
\end{Prop}

\begin{proof}
Again, we assume that the manifold $N$ is isometrically embedded in a $q$-dimensional
Euclidean vector space $\mathbb{R}^q$ and that the vector valued function \(u\) is a solution of 
(\ref{evolution-u-Rq}) and the spinor \(\psi\) a solution of (\ref{evolution-psi-Rq}).

Depending on the point of view, the function $u$ and the spinor \(\psi\) both satisfy an elliptic 
and a parabolic partial differential equation. 
This allows us to apply the classical Schauder estimates for both elliptic and parabolic equations.
First of all, $u$ satisfies the elliptic system
\begin{align*}
\Delta u=&\sff_u(u',u')+P(\sff_u(u',\partial_s\cdot\psi),\psi)
+\epsilon P(\sff_u(u',\nabla\psi),\psi)\\
&-\epsilon P(\sff_u(u',\psi),\nabla\psi)
+\epsilon B_u(u',\psi,u',\psi)+\dot{u}
\end{align*}
with the Laplacian $\Delta$ on $S^1$.
Let us estimate the right hand side
\begin{align*}
|\Delta u|\leq& C(|u'|^2+|\psi|^2|u'|+|\psi||\p\psi||u'|+|\dot{u}|+|\psi|^2|u'|^2)\leq C,
\end{align*}
where we used the estimates derived before.
Using Schauder estimates for solutions to an elliptic partial differential equation, we find
\begin{align*}
|u(\cdot,t)|_{C^{1+\alpha}(S^1,\R^q)}
&\leq C\big(\sup_{t\in [0,T)} |\Delta u(\cdot,t)|_{L^\infty(S^1,\R^q)}+\sup_{t\in[0,T)}|u(\cdot,t)|_{L^\infty(S^1,\R^q)}\big)\\
 &\leq C,
\end{align*}
since \(u\) takes values in a compact region of \(\R^q\).
We can improve the regularity of the spinor \(\psi\) by the same method.
Remember that $\psi\in\Gamma(\Sigma S^1\otimes T\R^q)$ solves the elliptic equation
\begin{align*}
\Delta\psi=\frac{1}{\epsilon}&\p\psi-\frac{1}{\epsilon}\sff_u(u',\partial_s\cdot\psi)
+\sff_u(\dot{u},\psi)+\frac{\nabla\psi}{\partial t}
+2\sff_u(u',\nabla\psi)\\
&+(\nabla\sff_u)(u',\psi)
+\sff_u(\tau(u),\psi)).
\end{align*}
Again, we bound the right hand side with the help of the previous estimates
\begin{align*}
|\Delta\psi|&\leq C\big(|\nabla\psi|+|u'||\psi|
+|\dot{u}||\psi|+\big|\frac{\nabla\psi}{\partial t}\big|+|u'||\nabla\psi|
+|u'|^2|\psi|+|u''||\psi|\big)\\
&\leq C
\end{align*}
and apply Schauder estimates for elliptic equations
\begin{align*}
&|\psi(\cdot,t)|_{C^{1+\alpha}(S^1,\Sigma S^1\otimes T\R^q)}\\ 
&\leq C\big(\sup_{t\in [0,T)} |\Delta \psi(\cdot,t)|_{L^\infty(S^1,\Sigma S^1\otimes T\R^q)}+\sup_{t\in[0,T)}|\psi(\cdot,t)|_{L^\infty(S^1,\Sigma S^1\otimes T\R^q)}\big)
\leq C.
\end{align*}
After exploiting the elliptic nature of the evolution equations for \((\gamma_t,\psi_t)\),
we now take the parabolic point of view.
Note that $u$ is also a solution of the parabolic system
\begin{align*}
Lu=&\sff_u(u',u')+P(\sff_u(u',\partial_s\cdot\psi),\psi)
+\epsilon P(\sff_u(u',\nabla\psi),\psi)\\
&-\epsilon P(\sff_u(u',\psi),\nabla\psi)
+\epsilon B_u(u',\psi,u',\psi)
\end{align*}
with $L=\Delta-\frac{\partial}{\partial t}$ denoting the heat operator on $S^1$.
Utilizing the previous estimate, we can bound the right hand side by
\begin{align*}
|\sff_u(du,du)+P(\sff_u&(u',\partial_s\cdot\psi),\psi)
+\epsilon P(\sff_u(u',\nabla\psi),\psi)\\
&-\epsilon P(\sff_u(u',\psi),\nabla\psi)+\epsilon B_u(u',\psi,u',\psi)|_{C^\alpha(S^1,\R^q)}\leq C.
\end{align*}
Finally, we employ Schauder estimates for linear parabolic partial differential equations
\begin{align*}
|u(\cdot,t)|&_{C^{2+\alpha}(S^1,\mathbb{R}^q)}+|\frac{\partial u}{\partial t}(\cdot,t)|_{C^\alpha(S^1,\mathbb{R}^q)}\\
&\leq C\big(\sup_{t\in [0,T)} |Lu(\cdot,t)|_{C^\alpha(S^1,\mathbb{R}^q)}+\sup_{t\in[0,T)}|u(\cdot,t)|_{L^\infty(S^1,\mathbb{R}^q)}\big)
\leq C,
\end{align*}
which proves the statement concerning the regularity of $u$.\\
Again, taking the parabolic point of view, \(\psi\) satisfies
\begin{align*}
L\psi=&\frac{1}{\epsilon}\p\psi-\frac{1}{\epsilon}\sff_u(u',\partial_s\cdot\psi)
+\sff_u(\dot{u},\psi)+\sff_u(u',\nabla\psi)
+2\sff_u(u',\nabla\psi)\\
&+(\nabla\sff_u)(u',\psi)+\sff_u(\tau(u),\psi)
\end{align*}
with the heat operator \(L=\Delta-\frac{\nabla}{\partial t}\).
The right hand side is bounded in \(C^\alpha\) such that 
\begin{align*}
|\psi&(\cdot,t)|_{C^{2+\alpha}(S^1,\Sigma S^1\otimes T\R^q)}+|\frac{\nabla\psi}{\partial t}(\cdot,t)|_{C^\alpha(S^1,\Sigma S^1\otimes T\R^q)}\\
&\leq C\big(\sup_{t\in [0,T)} |L\psi(\cdot,t)|_{C^\alpha(S^1,\Sigma S^1\otimes T\R^q)}+\sup_{t\in[0,T)}|\psi(\cdot,t)|_{L^\infty(S^1,\Sigma S^1\otimes T\R^q)}\big)
\leq C,
\end{align*}
which establishes the regularity of the spinor \(\psi\).
\end{proof}
Based on the estimates deduced so far, the uniqueness and stability result and the classical Schauder estimates,
we can now establish the long-time existence of the evolution equations.
\begin{Prop}[Long-time Existence]
Let 
\(
(\gamma_t,\psi_t)\in\chi(S^1\times[0,T),N)
\)
be a solution of (\ref{evolution-gamma-rescaled}) and (\ref{evolution-psi-rescaled})
and assume that \(N\) is compact.
Then for any initial data $(\gamma_0,\psi_0)\in C^{2+\alpha}(S^1,N)\times C^{2+\alpha}(S^1,\Sigma S^1\otimes\gamma_0^{-1}TN)$, 
there exists a unique 
\[
(\gamma_t,\psi_t)\in C^\infty(S^1\times(0,\infty),N)\times C^\infty(S^1\times(0,\infty),\Sigma S^1\otimes\gamma_t^{-1}TN)
\]
such that
\begin{align}
\label{curve-longtime-gamma}
&\begin{cases}
  \frac{\partial\gamma_t}{\partial t}=\tau(\gamma_t)-\cR(\gamma_t,\psi_t)-\epsilon\cR_c(\gamma_t,\psi_t),  & (s,t)\in S^1 \times (0,\infty),\\
  \gamma(s,0)=\gamma_0(s), &
\end{cases}\\
&\begin{cases}
\label{curve-longtime-psi}
  \frac{\tilde{\nabla}\psi_t}{\partial t}=\tilde{\Delta}\psi_t-\frac{1}{\epsilon}\D\psi_t,\qquad (s,t)\in S^1 \times (0,\infty),\\
  \psi(s,0)=\psi_0(s) & 
\end{cases}
\end{align}
holds.
\end{Prop}

\begin{proof}
The existence of a smooth short-time solution of the evolution equations is guaranteed by Theorem 3.22 in \cite{phd} for a time interval \(0\leq t\leq T_{max}\). 
We now demonstrate that, for $N$ being compact, the regularized Dirac-geodesic heat flow cannot blow up and will exist for all \(t\in[0,\infty)\).
We set
\[
T_0=\sup\{t\in[0,\infty)\mid(\ref{evolution-gamma-rescaled}),
(\ref{evolution-psi-rescaled})\text{ have a solution in } S^1\times[0,t)\}
\]
and show that $T_0=\infty$. Let us assume the opposite case.
We choose a sequence of numbers $\{t_i\}\subset [0,T_0)$ such that $t_i\rightarrow T_0$ as $i\rightarrow\infty$ 
and set $0<\alpha<\alpha'<1$. The embeddings
$C^{k+\alpha'}(S^1,N)\hookrightarrow C^{k+\alpha}(S^1,N)$ 
and in addition $C^{k+\alpha'}(S^1,\Sigma S^1\otimes\gamma^{-1}TN)\hookrightarrow C^{k+\alpha}(S^1,\Sigma S^1\otimes\gamma^{-1}TN)$ are compact, 
since \(S^1\) is compact.
By Proposition \ref{curve-prop-schauder}, the sequences
\[
\{\gamma(\cdot,t_i),\psi(\cdot,t_i)\}~\text{and}~\{\partial_t\gamma(\cdot,t_i),\tilde{\nabla}_t\psi(\cdot,t)\}
\]
are bounded in \(C^{2+\alpha'}(S^1,N)\times C^{2+\alpha'}(S^1,\Sigma S^1\otimes\gamma^{-1}TN)\) and 
\(C^{\alpha'}(S^1,N)\times C^{\alpha'}(S^1,\Sigma S^1\otimes\gamma^{-1}TN)\).
Hence, there exists a subsequence $\{t_{i_k}\}$ of $\{t_i\}$ with
\begin{align*}
&(\gamma(\cdot,T_0),\psi(\cdot,T_0))\in C^{2+\alpha'}(S^1,N)\times C^{2+\alpha'}(S^1,\Sigma S^1\otimes\gamma^{-1}TN),\\
&(\partial_t\gamma(\cdot,T_0),\tilde{\nabla}_t\psi(\cdot,T_0))\in C^{\alpha'}(S^1,N)\times C^{\alpha'}(S^1,\Sigma S^1\otimes\gamma^{-1}TN)
 \end{align*}
such that
\[
\{\gamma(\cdot,t_{i_k}),\psi(\cdot,t_{i_k})\}~\text{and}~\{\partial_t\gamma(\cdot,t_{i_k}),\tilde{\nabla}_t\psi(\cdot,t_{i_k})\}
\]
converge uniformly to $(\gamma(\cdot,T_0),\psi(\cdot,T_0))$ and $(\partial_t\gamma(\cdot,T_0),\tilde{\nabla}_t\psi(\cdot,T_0))$ respectively,
as $t_{i_k}\rightarrow T_0$. This is true for each \(t_{i_k}\) 
\begin{align*}
\frac{\partial\gamma}{\partial t}(\cdot,t_{i_k})=&(\tau(\gamma)-\cR(\gamma,\psi)-\epsilon\cR_c(\gamma,\psi))(\cdot,t_{i_k}), \\
\frac{\tilde{\nabla}\psi}{\partial t}(\cdot,t_{i_k})=&(\tilde{\Delta}\psi-\frac{1}{\epsilon}\D\psi)(\cdot,t_{i_k})
\end{align*}
and, consequently, also at \(T_0\).
Hence, the equations (\ref{curve-longtime-gamma}), (\ref{curve-longtime-psi}) have a solution in \(S^1\times[0,T_0]\).
We can now again apply the short-time existence Theorem 3.22 from \cite{phd} with initial values \((\gamma(\cdot,T_0),\psi(\cdot,T_0))\).
For \(\delta>0\) we then get a solution
\begin{eqnarray*}
&&\begin{cases}
  \frac{\partial\gamma_t}{\partial t}=\tau(\gamma_t)-\cR(\gamma_t,\psi_t)-\epsilon\cR_c(\gamma_t,\psi_t),  & (s,t)\in S^1 \times (T_0,T_0+\delta),\\
  \gamma(s,T_0)=\gamma_0(s), & 
\end{cases}\\
&&\begin{cases}
  \frac{\tilde{\nabla}\psi_t}{\partial t}(x)=\tilde{\Delta}\psi_t(x)-\frac{1}{\epsilon}\D\psi_t(x),\qquad (s,t)\in S^1 \times (T_0,T_0+\delta),\\
  \psi(s,T_0)=\psi_0(s) & 
\end{cases}
\end{eqnarray*}
in 
\begin{align*}
&\gamma\in C^{2+\alpha,1+\alpha/2}(S^1\times[T_0,T_0+\delta),N),\\
&\psi\in C^{2+\alpha,1+\alpha/2}(S^1\times[T_0,T_0+\delta),\Sigma S^1\otimes\gamma^{-1}TN).
\end{align*}
We realize that both solutions coincide on \(S^1\times\{T_0\}\) and for this reason,
we can glue them to a solution with existence interval \([0,T_0+\delta)\).
By standard regularity theory we find that this solution is smooth, for more details see Theorem 3.24 in \cite{phd}.
As a matter of fact, the system (\ref{curve-longtime-gamma}), (\ref{curve-longtime-psi}) has a smooth solution in \([0,T_0+\delta)\)
contradicting the definition of \(T_0\). The uniqueness of \((\gamma_t,\psi_t)\) then
follows from Proposition \ref{theorem-curve-uniqueness}.
\end{proof}

\section{Convergence}
In this section, we want to discuss under which assumptions and in which sense the evolution equations
for regularized Dirac-geodesics converge as \(t\to\infty\).
To this end, it is necessary to improve the estimates derived in the previous sections.
In particular, we need uniform estimates that do not depend on \(t\).
First of all, we will see that the issue of convergence depends crucially on
the norm of the spinor \(\psi_t\).
\begin{Prop}
Let 
\(
(\gamma_t,\psi_t)\in \chi(S^1\times[0,\infty),N)
\)
be a solution of (\ref{evolution-gamma-rescaled}) and (\ref{evolution-psi-rescaled}).
If we find a uniform bound on the spinor $\psi_t$, namely
\[
|\psi_t|^2\leq C,
\]
then we get a uniform bound on
\begin{equation}
|\gamma'_t|^2+\epsilon|\tilde{\nabla}\psi_t|^2\leq C
\end{equation}
for all \(t\in[0,\infty)\).
The constant $C$ depends on $N,\epsilon,\psi_0,\gamma'_0\) and \(\tilde{\nabla}\psi_0\).
\end{Prop}

\begin{proof}
By assumption, we have a uniform bound on $|\psi_t|^2$.
To derive the uniform bound on $|\gamma'_t|^2$ and $|\tilde{\nabla}\psi_t|^2$, we go back into the 
proof of Proposition \ref{curve-satz-evolution-F}. By the bound on $|\psi_t|^2$, we find that the quantity
$F_t:=\frac{1}{2}(|\gamma'_t|^2+\epsilon|\tilde{\nabla}\psi_t|^2)$ satisfies 
\[
\frac{\partial F_t}{\partial t}\leq\Delta F_t+CF_t.
\]
From the inequality \(E_\epsilon(\gamma_t,\psi_t)\leq E_\epsilon(\gamma_0,\psi_0)\), the bound on \(\psi_t\) and Young's inequality, 
we deduce
\[
\int_{S^1}(|\gamma'_t|^2+\epsilon|\tilde{\nabla}\psi_t|^2)ds\leq C.
\]
Applying (\ref{maximum-principle-l2}), the estimates on $|\gamma'_t|^2$ and $|\tilde{\nabla}\psi_t|^2$ follow.
\end{proof}
By the estimate just derived, we are now able to bound the \(t\) derivatives of \(\gamma_t\) and \(\psi_t\) uniformly.
\begin{Lem}
Let 
\(
(\gamma_t,\psi_t)\in \chi(S^1\times[0,\infty),N)
\)
be a solution of (\ref{evolution-gamma-rescaled}) and (\ref{evolution-psi-rescaled}).
If we can control the norm of $\psi_t$ uniformly, then we find
\begin{equation}
\big|\frac{\partial\gamma_t}{\partial t}\big|^2+\big|\frac{\tilde{\nabla}\psi_t}{\partial t}\big|^2\leq C
\end{equation}
for all \(t\in[0,\infty)\).
The constant $C$ depends on $N,\epsilon,\psi_0,\gamma'_0$ and $\tilde{\nabla}\psi_0$.
\end{Lem}
\begin{proof}
By Proposition \ref{curve-satz-evolution-G}, the quantity 
$G_t:=\frac{1}{2}\left(\big|\frac{\partial\gamma_t}{\partial t}\big|^2+\big|\frac{\tilde{\nabla}\psi_t}{\partial t}\big|^2\right)$
satisfies the inequality
\[
\frac{\partial G_t}{\partial t}\leq\Delta G_t+ Z(t)G_t.
\]
Applying the bounds on $|\psi_t|^2,|\gamma'_t|^2$ and $|\tilde{\nabla}\psi_t|^2$, we
find that $Z(t)$ is uniformly bounded such that \(G_t\) satisfies
\[
\frac{\partial G_t}{\partial t}\leq\Delta G_t+ CG_t.
\]
Integrating over $S^1$ and with respect to $t$ yields
\begin{align*}
\int_{S^1}\left(\big|\frac{\partial\gamma_t}{\partial t}\big|^2+\big|\frac{\tilde{\nabla}\psi_t}{\partial t}\big|^2\right)ds
\leq& C\int_0^\infty\int_{S^1}\left(\big|\frac{\partial\gamma_t}{\partial t}\big|^2+\big|\frac{\tilde{\nabla}\psi_t}{\partial t}\big|^2\right)dsdt
\\
&+\int_{S^1}\left(\big|\frac{\partial\gamma_t}{\partial t}\big|_{t=0}^2+\big|\frac{\tilde{\nabla}\psi_t}{\partial t}\big|_{t=0}^2\right)ds\\
\leq& E_\epsilon(\gamma_0,\psi_0)+C \\
\leq& C.
\end{align*}
The assertion follows from applying (\ref{maximum-principle-l2}) again.
\end{proof}

We realize that the convergence of the evolution equations depends crucially on the norm of \(\psi_t\). 
One way of controlling the norm of $\psi_t$ is to choose the parameter $\epsilon$
large enough such that the parabolic nature of the evolution equation for $\psi_t$
dominates, which basically means that the second order term is sufficiently
large to control the first order one.
\begin{Prop}
Let 
\(
(\gamma_t,\psi_t)\in \chi(S^1\times[0,\infty),N)
\)
be a solution of (\ref{evolution-gamma-rescaled}) and (\ref{evolution-psi-rescaled}).
For $\epsilon\geq 1$ we get a uniform bound of $|\psi_t|^2$ for all \(t\in[0,\infty)\).
\end{Prop}
\begin{proof}
Using the evolution equation (\ref{evolution-psi-rescaled}), we calculate
\begin{align*}
\frac{\partial}{\partial t}\frac{1}{2}\int_{S^1}|\psi_t|^2ds=&-\int_{S^1}|\tilde{\nabla}\psi_t|^2ds
+\frac{1}{\epsilon}\int_{S^1}\langle\psi_t,\D\psi_t\rangle ds \\
\leq&-\int_{S^1}|\tilde{\nabla}\psi_t|^2ds+\frac{1}{\epsilon}\int_{S^1}|\psi_t||\tilde{\nabla}\psi_t|ds\\
\leq&(\frac{1}{\epsilon}-1)\int_{S^1}|\tilde{\nabla}\psi_t|^2ds,
\end{align*}
where we used the Cauchy-Schwarz and the Poincaré inequality on $S^1$ in the last step.
After integration with respect to \(t\) we find for $\epsilon\geq 1$ 
\[
\int_{S^1}|\psi_t|^2ds\leq \int_{S^1}|\psi_0|^2ds.
\]
We have already seen that $|\psi_t|^2$ satisfies the pointwise equation
\[
\frac{\partial}{\partial t}\frac{1}{2}|\psi_t|^2\leq\Delta\frac{1}{2}|\psi_t|^2+\frac{1}{4\epsilon^2}|\psi_t|^2,
\]
and by (\ref{maximum-principle-l2}) we get a uniform bound on $|\psi_t|^2$.
\end{proof}

\begin{Bem}
Of course, it would be much nicer if one could bound the norm of \(\psi_t\)
while keeping the regularizing parameter \(\epsilon\) small. Unfortunately, this does not seem to be possible.
From a spectral point of view \(\epsilon\geq 1\) implies that the operator
\(
L=\epsilon\tilde{\Delta}-\D
\)
is positive.
\end{Bem}

After having derived pointwise uniform bounds, the regularity of $(\gamma_t,\psi_t)$
can be improved by applying Schauder estimates again.

The existence of a convergence subsequence of the evolution equations can now be derived by standard methods.
\begin{Lem}[Convergence]
Let
\(
(\gamma_t,\psi_t)\in \chi(S^1\times[0,\infty),N)
\)
be a solution of (\ref{evolution-gamma-rescaled}) and (\ref{evolution-psi-rescaled}).
If $\epsilon\geq 1$ then $(\gamma_t,\psi_t)$ subconverges to a regularized Dirac-geodesic $(\gamma_\infty,\psi_\infty)$ in
\(C^2(S^1,N)\times C^2(S^1,\Sigma S^1\otimes\gamma_t^{-1}TN)\).
\end{Lem}
\begin{proof}
First of all, we improve the regularity of our estimates with the help of Schauder theory, as in Proposition \ref{curve-prop-schauder} and find
\begin{align*}
&\sup_{t\in [0,\infty)}\left(|\gamma(\cdot,t)|_{C^{2+\alpha}(S^1,N)}+\big|\frac{\partial\gamma}{\partial t}(\cdot,t)\big|_{C^\alpha(S^1,N)}\right)\leq C,\\
&\sup_{t\in [0,\infty)}\left(|\psi(\cdot,t)|_{C^{2+\alpha}(S^1,\Sigma S^1\otimes\gamma_t^{-1}TN)}+\big|\frac{\tilde{\nabla}\psi}{\partial t}(\cdot,t)\big|_{C^{\alpha}(S^1,\Sigma S^1\otimes\gamma_t^{-1}TN)}\right)\leq C,
\end{align*}
but now the constants \(C\) do not depend on \(t\).
In addition, we have the estimate from the inequality for the energy \(E_\epsilon(\gamma,\psi)\)
\[
\int_0^{\infty}\int_{S^1}\left(\big|\frac{\partial\gamma_t}{\partial t}\big|^2
+\big|\frac{\tilde{\nabla}\psi_t}{\partial t}\big|\right)dsdt\leq C.
\]
Hence, there exists a subsequence $t_k$ such that as \(k\to\infty\) we have
\[
\big|\frac{\partial\gamma(\cdot,t_k)}{\partial t}\big|^2_{L^2(S^1\times [0,\infty))}\to 0,\qquad \big|\frac{\tilde{\nabla}\psi(\cdot,t_k)}{\partial t}\big|^2_{L^2(S^1\times [0,\infty))}\to 0.
\]
Using Schauder estimates again,
\begin{align*}
\sup_{t\in [0,\infty)}\left(|\gamma(\cdot,t_k)|_{C^{2+\alpha}(S^1,N)}+\big|\frac{\partial\gamma}{\partial t}(\cdot,t_k)\big|_{C^\alpha(S^1,N)}\right)\leq C,\\
\sup_{t\in [0,\infty)}\left(|\psi(\cdot,t_k)|_{C^{2+\alpha}(S^1,\Sigma S^1\otimes\gamma_t^{-1}TN)}+\big|\frac{\tilde{\nabla}\psi}{\partial t}(\cdot,t_k)\big|_{C^{\alpha}(S^1,\Sigma S^1\otimes\gamma_t^{-1}TN)}\right)\leq C,
\end{align*}
it follows from the Theorem of Arzela Ascoli that there exists a convergent subsequence,
which is also denoted by $t_k$, such that the pair $(\gamma_{t_k},\psi_{t_k})$ converges in 
the space $C^2(S^1,N)\times C^2(S^1,\Sigma S^1\otimes\gamma_t^{-1}TN)$ to
a limiting map $(\gamma_\infty,\psi_\infty)\).
Since $(\gamma_t,\psi_t)$ is smooth in $t$, we find that $(\gamma_\infty,\psi_\infty)$ is
homotopic to $(\gamma_0,\psi_0)$.
\end{proof}
The smoothness of the limiting map \((\gamma_\infty,\psi_\infty)\) follows from elliptic estimates
as in Proposition \ref{curve-prop-schauder}.

Let us make several remarks on what we have achieved so far.

\begin{Bem}
It would be desirable to obtain further properties of the limiting map \((\gamma_\infty,\psi_\infty)\).
In particular, it would be nice to ensure that the limit map \((\gamma_\infty,\psi_\infty)\) is not trivial.
In the case of the heat flow for geodesics it can happen that the initial curve \(\gamma_0\) is shrunk to
a point under the evolution.
Of course, one would like to exclude this possibility here.
\end{Bem}

\begin{Bem}
In the case of the harmonic map heat flow and the assumption $K^N\leq 0$, it
is known that the limit $k\to\infty$ is independent of the chosen subsequence.
This result is known as Hartmann's theorem \cite{MR0214004} and makes use
of the fact that the second variation of the energy functional is positive.
In the case of regularized Dirac-geodesics, the second variation of the energy
functional \(E_\epsilon(\gamma,\psi)\) is not positive (see \cite{phd}, Proposition 2.2) and we cannot 
derive an analogue of Hartmann's theorem.
\end{Bem}

\begin{Bem}
A question connected to the last remark is if the gradient flow converges itself or if it 
just posseses a convergent subsequence. Although this issue seems simple at first glance
it turns out to be rather subtle.
For the geodesic heat flow this question was recently resumed in \cite{choi-parker}.
In that reference the authors give examples for the geodesic heat flow 
in which the energy functional is not convex and the heat flow does not converge as \(t\to\infty\),
although having a convergent subsequence.
\end{Bem}

\section{Removing the Regularization}
By the considerations we did so far, we have obtained a \emph{regularized Dirac-geodesic}.
Namely, we constructed a smooth solution to the coupled problem
\begin{align}
\tau(\gamma_\infty)=&\cR(\gamma_\infty,\psi_\infty)+\epsilon\cR_c(\gamma_\infty,\psi_\infty),\\
\epsilon\tilde{\Delta}\psi_\infty=&\D\psi_\infty.
\end{align}
To obtain a \emph{Dirac-geodesic}, we have to remove the regularization and let \(\epsilon\to 0\).
Aiming in this direction, the next Lemma shows how the critical points of \(E(\gamma,\psi)\) and \(E_\epsilon(\gamma,\psi)\)
are related to each other. Note that the spectrum of the twisted Dirac operator \(\D\) is real and discrete.

\begin{Lem}[Critical points of $E_\epsilon(\gamma,\psi)$ and $E(\gamma,\psi)$]
For $\epsilon\neq-\frac{1}{\lambda}$, where $\lambda$ is an eigenvalue of
the twisted Dirac-operator $\D$, the regularized functional  $E_\epsilon(\gamma,\psi)$ 
has the same critical points as $E(\gamma,\psi)$.
\end{Lem}
\begin{proof}
The critical points of the functional $E(\gamma,\psi)$ are given by
\begin{equation*}
\tau(\gamma)=\cR(\gamma,\psi),\qquad \D\psi=0,
\end{equation*}
whereas the critical points of the regularized functional $E_\epsilon(\gamma,\psi)$ are given by
\begin{equation*}
\tau(\gamma)=\cR(\gamma,\psi)+\epsilon\cR_c(\gamma,\psi),\qquad \D\psi+\epsilon\D^2\psi=0.
\end{equation*}
It is clear, that if $(\gamma,\psi)$ is a Dirac-geodesic, then it is also a regularized Dirac-geodesic.\\
The other direction is slightly more subtle.
Assume that $(\gamma,\psi)$ is a regularized Dirac-geodesic. Using the equation for $\psi$
and integrating over $S^1$, we obtain
\[
\int_{S^1}\langle\psi,\D\psi\rangle ds+\epsilon\int_{S^1}|\D\psi|^2ds=0.
\]
For this equation to hold, either $\psi$ must be trivial or $\psi+\epsilon\D\psi=0$.
But we have chosen $\epsilon$ in such a way that the second possibility is excluded.
We conclude that \(\D\psi=0\). Hence, \(\cR_c(\gamma,\psi)=0\)
and therefore the pair $(\gamma,\psi)$ is a Dirac-geodesic.
\end{proof}
Thus, we have to ensure that we do not lose regularity in the limit \(\epsilon\to 0\).
It is easy to see that all estimates that were derived when studying the evolution equations for \((\gamma,\psi)\)
do not survive the limit \(\epsilon\to 0\). Therefore the question we have to ask is, 
which estimates hold independently of \(\epsilon\). \\
To this end, let us compute
\begin{equation}
\frac{\partial}{\partial s}\frac{1}{2}|\gamma'_\infty|^2=\langle\gamma'_\infty,\tau(\gamma_\infty)\rangle
=\langle\gamma'_\infty,\cR(\gamma_\infty,\psi_\infty)+\epsilon\cR_c(\gamma_\infty,\psi_\infty)\rangle=0.
\end{equation}
and thus the norm of \(\gamma'_\infty\) is constant, independently of \(\epsilon\). \\
In the following we again assume that \(N\subset\R^q\).
After taking the limit \(\epsilon\to 0\), the pair \((u,\psi)\) solves the following set of equations,
\begin{align}
\label{u-infty} u''=&\sff(u',u')+P(\sff(\partial_s\cdot\psi,u'),\psi),\\
\label{psi-infty}\p\psi_\infty=&\sff(\partial_s\cdot\psi,u').
\end{align}
The right hand side of (\ref{u-infty}) is in \(L^\infty\) by the bound on \(\gamma'_\infty\) and the fact that the spinor \(\psi_\infty\) can be interpreted as a function on \(S^1\).
Hence, by the usual Schauder theory we get that \(u\in C^{1+\alpha}(S^1,\R^q)\). Differentiating (\ref{psi-infty}) and estimating the right hand side, we also find that
\(\psi\in C^{1+\alpha}(S^1,\Sigma S^1\otimes T\R^q)\). Then a standard bootstrap argument yields that the pair \((u,\psi)\) is smooth and consequently,
we cannot lose regularity in the limit \(\epsilon\to 0\).

Nevertheless the properties of the pair \((\gamma_\infty,\psi_\infty)\) after \(\epsilon\to 0\) have to be
further investigated.

\appendix
\section{}
The following Lemma combines the pointwise maximum principle 
with an integral norm. It can be thought of as a simple version of 
Moser's parabolic Harnack inequality.
\begin{Lem}
\label{maximum-principle-l2}
Assume that \((M,h)\) is a compact Riemannian manifold. If a function \(u(s,t)\geq 0\) satisfies
\[
\frac{\partial u}{\partial t}\leq \Delta u+Cu,
\]
and if in addition we have the bound
\[
U(t)=\int_Mu(s,t)dM\leq U_0,
\]
then there exists a uniform bound on 
\[
u(s,t)\leq e^CKU_0 
\]
with the constant \(K\) depending on \(M\).
\end{Lem}
\begin{proof}
A proof can for example be found in \cite{MR2744149}, p.\ 284.
\end{proof}

\emph{Acknowledgements:}
The author would like to thank the ``IMPRS for Geometric Analysis, Gravitation and String Theory'' 
for financial support.

\bibliographystyle{amsplain}
\bibliography{mybib}
\end{document}